\documentclass[11pt]{amsart}

\usepackage{amsfonts, amstext, amsmath}
\usepackage{graphicx, color, epsfig}

\newcommand{\D}{{\partial}} 

\newcommand{\vol}{\mathop{\rm vol}\nolimits}

\newcommand{\TT}{{\mathbb{T}}}
\newcommand{\CN}{{\mathcal{N}}}
\newcommand{\CT}{{\mathcal{T}}}

\theoremstyle{plain}
\newtheorem{theorem}{Theorem}[section]
\newtheorem{corollary}[theorem]{Corollary}
\newtheorem{lemma}[theorem]{Lemma}
\newtheorem{prop}[theorem]{Proposition}

\newtheorem*{namedtheorem}{\theoremname}
\newcommand{\theoremname}{testing}
\newenvironment{named}[1]{\renewcommand{\theoremname}{#1}\begin{namedtheorem}}{\end{namedtheorem}}

\theoremstyle{definition}
\newtheorem{define}[theorem]{Definition}
\newtheorem{remark}[theorem]{Remark}

\title[Multiply twisted knots that are Sefert fibered or toroidal]{On multiply
	twisted knots that are \\Seifert fibered or toroidal}

\author{Jessica S. Purcell}

\address{Jessica S. Purcell, Department of Mathematics, Brigham Young
University, Provo, UT 84602}

\email{jpurcell@math.byu.edu}

\begin{document}
\bibliographystyle{hamsplain}

\begin{abstract}
	We consider knots whose diagrams have a high amount of twisting of
multiple strands.  By encircling twists on multiple strands with
unknotted curves, we obtain a link called a generalized augmented
link.  Dehn filling this link gives the original knot.  We classify
those generalized augmented links that are Seifert fibered, and give a
torus decomposition for those that are toroidal.  In particular, we
find that each component of the torus decomposition is either
``trivial'', in some sense, or homeomorphic to the complement of a
generalized augmented link.  We show this structure persists under
high Dehn filling, giving results on the torus decomposition of knots
with generalized twist regions and a high amount of twisting.  As an
application, we give lower bounds on the Gromov norms of these knot
complements and of generalized augmented links.  

\end{abstract}

\maketitle


\section{Introduction}
\label{sec:intro}
This paper continues a program to understand the geometry of knot and
link complements in $S^3$, given only a diagram of the knot or link.
Each knot complement decomposes uniquely along incompressible tori
into hyperbolic and Seifert fibered pieces, by work of Jaco--Shalen
\cite{jaco-shalen} and Johannson \cite{johannson}.  By Mostow--Prasad
rigidity, the metric on the hyperbolic pieces is unique.  Thus this
geometric information on the complement is completely determined by a
diagram of the knot.  However, reading geometric information off of a
diagram seems to be difficult.

In recent years, techniques have been developed to relate geometric
properties to a diagram for classes of knots and links admitting
particular types of diagrams, such as alternating
\cite{lackenby:alt-volume}, and highly twisted knots and links
\cite{purcell:cusps, purcell:volume, fkp}.  However, many links of
interest to knot theorists and hyperbolic geometers do not admit these
types of diagrams.  These include Berge knots \cite{berge, baker:I,
baker:II}, twisted torus knots and Lorenz knots \cite{birman-kofman},
which contain many of the smallest volume hyperbolic knots
\cite{champanerkar-kofman}.  These knots admit diagrams that are
highly non-alternating, that have few twists per twist region, but
contain regions where multiple strands of the diagram twist around
each other some number of times.  The ideas of this paper and a
companion paper \cite{purcell:aug} grew out of a desire to understand
geometric properties of these ``multiply twisted'' knots and links,
given only a diagram.  The results here give a first step towards such
an understanding.

In this paper, we consider a class of multiply twisted knots in $S^3$,
described below, and classify those which are \emph{not} hyperbolic.
We completely determine those which are Seifert fibered, and describe
the unique torus decomposition (the JSJ decomposition) for those which
are toroidal.  We obtain our results by \emph{augmenting} diagrams of
the multiply twisted knots, that is, encircling regions of the diagram
where multiple strands twist about each other by a simple closed
curve, called a \emph{crossing circle}. This generalizes a
construction of Adams \cite{adams:aug}.  Since all knots are obtained
by Dehn filling some generalized augmented link, the geometric
properties of these links are also interesting.  In
\cite{purcell:aug}, we describe geometric properties of those
generalized augmented links which are known to be hyperbolic.
Combining these results leads to the JSJ decomposition of certain
multiply twisted knots.

To state our results carefully, we need some definitions. 

\subsection{Generalized augmented links}

First, although our main results relate to knots in $S^3$, in fact
many results in this paper apply more generally to knots and links in
a 3--manifold $M$.  We will assume $M$ is compact, with (possibly
empty) boundary consisting of tori, and $M$ admits an orientation
reversing involution fixing a surface $S$.  For example, $S^3$ is such
a manifold, taking $S$ to be a separating 2--sphere.  A solid torus
$V$ is another example, with $S$ a M\"obius band or annulus.

Let $K$ be a knot or link in $M$ that can be ambient isotoped into a
neighborhood of $S$.  We define a \emph{diagram} of the knot or link
$K$ with respect to the surface $S$ to be a projection of $K$ to $S$
yielding a 4--valent graph on $S$ with over--under crossing
information at each vertex.

Given a diagram, we may define twisting, twist regions, and
generalized twist regions exactly as in \cite{purcell:aug}, whether or
not our link is in $S^3$.  We review the definitions briefly here.
More precise statements are found in \cite{purcell:aug}.

\begin{figure}
	(a)
	\includegraphics{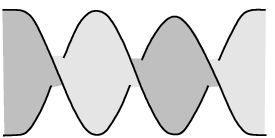}
	\hspace{.2in}
	(b)
	\includegraphics{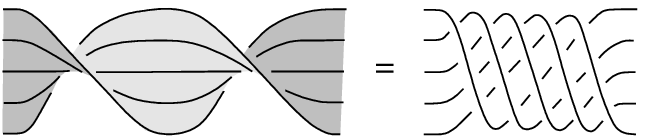}
	\caption{(a) A twist region.  (b) A generalized twist region.
	Multiple strands lie on the twisted ribbon surface.}
\label{fig:twist}
\end{figure}

A \emph{twist region} of a diagram is a region in which two strands
twist about each other maximally, as in Figure \ref{fig:twist}(a).
Note that the two strands bound a ``ribbon surface'' between them.  A
\emph{generalized twist region} is a region of a diagram in which
multiple strands twist about each other maximally, as in Figure
\ref{fig:twist}(b).  Note that all strands lie on the ribbon surface
bounded between the outermost strands.  A \emph{half--twist} of a
generalized twist region consists of a single crossing of the two
outermost strands, which flips the ribbon surface over once.  Figure
\ref{fig:twist}(b) shows a single full--twist, or two half--twists of
five strands.

We may group all crossings of a diagram into generalized twist
regions, so that each crossing is contained in exactly one generalized
twist region.  This is called a \emph{maximal twist region selection},
and is not necessarily unique.  For example, in Figure
\ref{fig:twist}(b), we could group all crossings into one generalized
twist region on five strands, or group each crossing into its own
twist region on two strands.

\begin{figure}
\begin{tabular}{ccccc}
	(a) & 
	\includegraphics{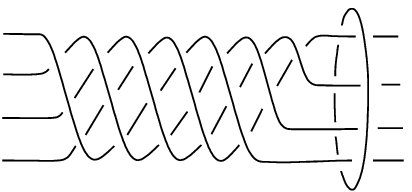} & \hspace{.2in} & (b) &
	\includegraphics{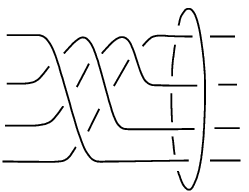}
\end{tabular}
\caption{(a) Encircle each twist region with a crossing circle.  (b)
	Link $L$ given by removing full--twists from the diagram.}
\label{fig:cross-cir}
\end{figure}

Given a maximal twist region selection, at each generalized twist
region insert a \emph{crossing circle}, i.e. a simple closed curve
$C_i$ encircling the strands of the generalized twist region, bounding
a disk $D_i$ perpendicular to the projection plane of the diagram.
The complement of the resulting link is homeomorphic to the complement
of the link whose diagram is obtained by removing all full--twists
from each generalized twist region.  This is illustrated in Figure
\ref{fig:cross-cir}.  The resulting link, with crossing circles added
and all full--twists removed, is defined to be a \emph{generalized
augmented link}.  We will always assume such a link contains at least
one crossing circle, to avoid trivial cases.  

Note that if $L$ is a generalized augmented link, obtained by
augmenting a knot $K$ in $M$, then $M-K$ is obtained from $M-L$ by
Dehn filling.  Let $\CN(C_i)$ denote a small embedded tubular
neighborhood of $C_i$ in $M$.  Let $\mu_i$ be (the isotopy class of)
the meridian of $\CN(C_i)$ (i.e. $\mu_i$ bounds a disk in $M$), and
let $\lambda_i = \D D_i$ be the longitude.  Suppose $n_i$ full--twists
were removed at $C_i$ to go from the diagram of $K \cup (\cup\, C_j)$
to that of $L$.  Then Dehn filling along the slope
$\mu_i+n_i\lambda_i$ on $\D \CN(C_i)$, for each $i$, yields $M-K$.
See, for example, Rolfsen \cite{rolfsen} for a more complete
description of this process.

We refer to this type of Dehn filling as \emph{twisting} along the
disk $D_i$, or along $C_i$, or, when $D_i$ or $C_i$ are understood,
simply as \emph{twisting}.  Note any link in $S^3$ is obtained by
twisting some generalized augmented link.

Finally, we wish to use diagrams of knots that do not involve
unnecessary twisting.  That is, we wish them to be \emph{reduced} in
the sense of the following definition, which we will use to generalize
Lackenby's definition of twist reduced \cite{lackenby:alt-volume}, and
Menasco--Thistlethwaite's definition of standard
\cite{menasco-thistle}.

\begin{define}
A generalized augmented link with knot strands $K_j$ and crossing
circles $C_i$ is said to be \emph{reduced} if the following hold:
\begin{enumerate}
	\item (Minimality of twisting disks) The twisting disk $D_i$ with
	boundary $C_i$ intersects $\cup\, K_j$ in $m_i$ points, where
	$m_i\geq 2$, and $m_i$ is minimal over all disks in $M$ bounded by
	$C_i$.  That is, if $E_i$ is another disk embedded in $M$ with
	boundary $C_i$, disjoint from all other crossing circles, then $|E_i
	\cap (\cup K_j)| \geq m_i$.
	\item (No redundant twisting)  There is no annulus embedded in $M-L$
	with one boundary component isotopic to $\D D_j$ on $\D \CN(C_j)$ and
	the other isotopic to $\D D_i$ on $\D \CN(C_i)$, $i\neq j$.
	\item (No trivial twisting)  There is no annulus embedded in $M-L$
	with one boundary component isotopic to $\D D_j$ on $\D \CN(C_j)$ and
	the other boundary component on $\partial M$.  
\end{enumerate}
\label{def:reduced}
\end{define}

\begin{figure}
\input{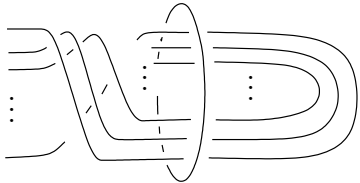}
\caption{Untwist at a homotopically trivial crossing circle.}
\label{fig:nugatory}
\end{figure}

\begin{figure}
	\input{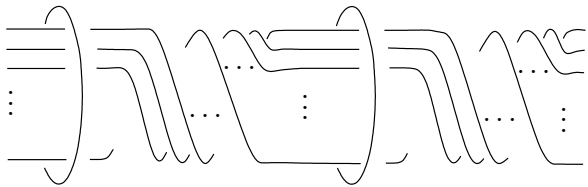}
\caption{Concatenate twists when crossing circles bound an annulus.}
\label{fig:redundant}
\end{figure}

These conditions allow us to rule out unnecessary crossing circles and
generalized twist regions.  For example, condition (1) prohibits
``nugatory'' twist regions, such as those shown in Figure
\ref{fig:nugatory}.  Condition (2) rules out redundant generalized
twist regions and their associated crossing circles, such as shown in
Figure \ref{fig:redundant}, where two crossing circles encircle the
same generalized twist region.

\subsection{Results}

We are now ready to state the main results of this paper.

In Section \ref{sec:sf}, we classify all reduced augmented links which
are Seifert fibered.  This is the content of Theorem \ref{thm:no-sf}.
As a consequence, we obtain the following result.

\begin{named}{Corollary \ref{cor:sf}}
Let $K$ be a knot in $S^3$ which has a diagram $D$ whose augmentation
is a Seifert fibered reduced augmented link.  Then $K$ is a $(2,q)$
torus knot.
\end{named}

In Section \ref{sec:satellite}, we describe a torus decomposition of
generalized augmented links.  In particular, in Theorem
\ref{thm:tori}, we show that components of such a torus decomposition
are either in some sense trivial, or atoroidal reduced augmented
links.  Applying this to links in $S^3$, using results on hyperbolic
generalized augmented links of \cite{purcell:aug}, we obtain a torus
decomposition for all links in $S^3$ obtained by high Dehn fillings of
toroidal generalized augmented links.  This is Theorem
\ref{thm:dehnfill}.  In particular, if at least 6 half--twists are
inserted when we twist along $C_i$, then the torus decomposition of
$M-K$ will agree with that of $M-L$, aside from ``trivial'' pieces.

We wish to apply these results to as many knots in $S^3$ as possible.
In Section \ref{sec:diagram}, we show that any knot in $S^3$ admits a
diagram whose augmentation is reduced.  Thus the above results will
apply at least to their augmentations.  This is Theorem
\ref{thm:reduced-knot}. 

In Section \ref{sec:knot}, we apply these results to knots in $S^3$,
to determine the JSJ decomposition of multiply twisted knots.  Our
main result is the following.

\begin{named}{Theorem \ref{thm:toroidal}}
Let $K$ be a knot in $S^3$ which is toroidal, with a twist--reduced
diagram and a maximal twist region selection with at least $6$
half--twists in each generalized twist region.  Let $L$ denote the
corresponding augmentation.  Then there exists a sublink $\hat{L}$ of
$L$, possibly containing fewer crossing circles, such that:
\begin{enumerate}
	\item The essential tori of the JSJ decomposition of $S^3-K$ are in
	one--to--one correspondence with those of $S^3-\hat{L}$.
	\item Corresponding components of the torus decompositions have the
	same geometric type, i.e. are hyperbolic or Seifert fibered.
	\item Essential tori of $S^3-\hat{L}$ and $S^3-K$ form a collection
	of nested tori, each bounding a solid torus in $S^3$ which contains
	$K$, and is fixed under a reflection of $S^3-L$.
\end{enumerate}
\end{named}

Putting this theorem with the results on hyperbolic geometry of
generalized augmented links in \cite{purcell:aug}, we obtain as an
application a lower bound on the Gromov norm of such knots.

\begin{named}{Theorem \ref{thm:gromov-norm-knots}}
Let $K$ be a knot in $S^3$ which is toroidal, with a twist--reduced
diagram at least $7$ half--twists in each generalized twist region.
Let $L$ denote the corresponding augmentation, and let $\hat{L}$
denote the sublink of Theorem \ref{thm:toroidal}.  Let $t$ denote the
number of crossing circles of $\hat{L}$.  Then the Gromov norm of
$S^3-K$ satisfies
$$\| [S^3-K] \| \geq 0.65721\,(t-1).$$
\end{named}

\subsection{Comments and additional questions}

The results of this paper give geometric information based purely on
diagrammatical properties of extensive classes of knots.  However,
because of the high amount of twisting required, for example in
Theorem \ref{thm:toroidal}, these classes still do not include
examples of many knots.  Considering knots which are not included
leads to two interesting remaining questions.

First, can the results of Theorem \ref{thm:toroidal} be sharpened to
require fewer half--twists?  We can construct examples of atoroidal
knots $K$ whose geometric type (hyperbolic or Seifert fibered) does
not agree with that of the corresponding reduced augmented link $L$.
However, in all these examples, at least one generalized twist region
contains fewer than $3$ half--twists.  In \cite{amm:twisting},
A{\"i}t-Nouh, Matignon, and Motegi, working on a related question,
show that when exactly one crossing circle is inserted into the
diagram of an unknot, and then the unknot is twisted, inserting at
least $4$ half--twists, the geometric type of the resulting knot
(Seifert fibered, toroidal, or hyperbolic) agrees with that of the
unknot union the crossing circle.  While these results do not apply to
generalized augmented links, the result requiring only $4$
half--twists is intriguing.

Secondly, given an arbitrary diagram of a knot, is there a way to
optimize the maximal twist region selection?  There are many ways to
choose a maximal twist region selection.  For example, Figure
\ref{fig:twist}(b) could either be seen as a single full--twist of $5$
strands, or as $20$ half--twists, each of $2$ strands.  To apply the
results of this paper, it seems we would want to select generalized
twist regions to maximize the number of half--twists in each twist
region.  Is there an algorithm that, given a diagram of a knot $K$,
produces a reduced diagram and a maximal twist region selection with
the highest number of half--twists per generalized twist region
possible for $K$?  Results along these lines would be interesting.

\subsection{Acknowledgements}
This research was supported in part by NSF grant DMS--0704359.  We
thank John Luecke for helpful conversations.

\section{Reflection}
\label{sec:reflect}
A generalized augmented link in a manifold $M$ admits a reflection
through a surface fixed pointwise, just as in the case of links in
$S^3$ \cite[Proposition 3.1]{purcell:slopes}.  This reflection is
necessary for many of the results that follow, and so we state it
first.

\begin{prop}
Let $M$ be a 3--manifold with torus boundary which admits an
orientation reversing involution fixing a surface $S$ pointwise.  Let
$K$ be a link in $M$ which may be isotoped to lie in a neighborhood of
$S$.  Finally, let $L$ be an augmentation of a diagram of $K$, with
crossing circles $\{C_1, \dots, C_n\}$ and knot strands $\{K_1, \dots,
K_m\}$.  Then $M-(\cup\, C_i)$ admits an orientation reversing
involution $\sigma$ which fixes a surface $P$ pointwise, and each
$K_j$ is embedded in $P$.  In particular, $M-L$ admits an orientation
reversing involution $\sigma$.
\label{prop:reflect}
\end{prop}

\begin{figure}
	\begin{center}
	\includegraphics{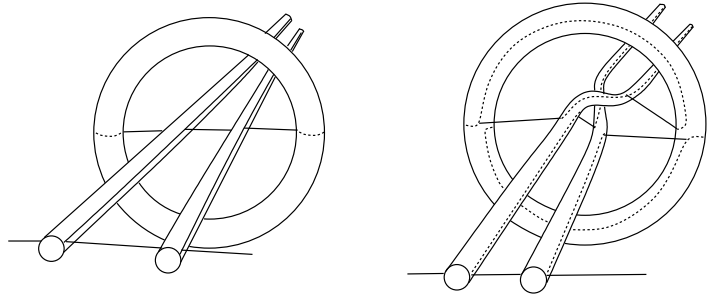}
	\end{center}
\caption{Left: $P\cap \D \CN(C_i)$ consists of two meridians when there
	is no half--twist.  Right: Under a half--twist, $P\cap \D \CN(C_i)$
	has boundary shown by the dotted lines. }
\label{fig:half-twist}
\end{figure}

\begin{proof}
Isotope crossing circles to be orthogonal to $S$, preserved by the
reflection of $M$ through $S$.  If there are no half--twists in the
diagram of $L$, then all components $K_j$ of $L$ are embedded in $S$.
Hence the reflection in $S$ preserves each $K_j$ as well as each
$C_j$, so $P=S$ and the involution is the restriction of the
involution of $M$ to $M-L$.

If there are half--twists in the diagram of $L$, then the reflection
of $M$ through the surface $S$ gives a new link $L'$ in which all the
directions of the crossings at each half--twist have been reversed.
Let $\tau$ be the homeomorphism of $M-(\cup\, C_j)$ which twists
exactly one full time in the opposite direction of these half--twists
at each corresponding crossing circle.  Applying $\tau$ changes the
diagram to one in which the crossings of half--twists have been
reversed again, hence to the diagram of $L$.  So $M-L'$ is
homeomorphic to $M-L$, and the orientation reversing involution
$\sigma$ of $M-(\cup\, C_j)$ is given by reflection of $M$ in $S$
followed by the homeomorphism $\tau$.

Finally, we describe the surface $P$ fixed pointwise by $\sigma$ in
the case of half--twists.  In this case, $P$ is equal to $S$ outside a
neighborhood of those crossing disks for which the corresponding
crossing circle $C_i$ bounds a half--twist.  Inside such a
neighborhood, the surface $P$ follows the ribbon surface of the
half--twist between the outermost knot strands.  Between $C_i$ and the
outermost knot strands, $P$ runs over the overcrossing, under the
undercrossing, and meets up with the surface $S$ on the opposite side
of the link.  Its boundary $P \cap \D \CN(C_i)$ runs twice along the
meridian of $\D \CN(C_i)$, once along the longitude, as in Figure
\ref{fig:half-twist}.
\end{proof}

In \cite[Lemma 3.1]{purcell:aug}, we showed that in this setting, the
slopes $\D D_i$ on $\D \CN(C_i)$ will meet the surface $P$ of
Proposition \ref{prop:reflect} exactly twice.  We will use this fact
here, and so we state it as a lemma.

\begin{lemma}
Let $L$ be a generalized augmented link in $M$, with reflective
surface $P$ from Proposition \ref{prop:reflect} and twisting disks
$D_i$.  Then for each $i$, $\D D_i$ meets $P$ exactly twice on $\D
\CN(C_i)$.\qed
\label{lemma:meet-twice}
\end{lemma}

\section{The Seifert fibered case}
\label{sec:sf}
In this section, we classify Seifert fibered augmented links.

\subsection{Incompressibility of surfaces}

Let $M-L$ be the complement of a generalized augmented link.  By
Proposition \ref{prop:reflect}, $M-L$ admits an involution $\sigma$
which fixes a surface $P$.  We show that $P$ and the (punctured)
twisting disks $D_i-(D_i \cap L) \subset M-L$ are incompressible.

\begin{lemma}
Let $N$ be an irreducible $3$--manifold with torus boundary components
which admits an orientation reversing involution $\sigma$ with fixes a
surface $P$.  Then $P$ is incompressible.
\label{lemma:P-incompress}
\end{lemma}

\begin{proof}
Suppose not.  Suppose $D$ is a compressing disk for $P$.  Then
$\partial D$ lies on $P$, so is fixed by $\sigma$, and $\sigma(D)$ is
a disk whose boundary agrees with that of $D$.  We first show we can
assume $D$ and $\sigma(D)$ are disjoint except on their boundaries.
If not, consider the intersections of $D$ and $\sigma(D)$.  Note these
intersections consist of closed curves on $D$ and on $\sigma(D)$ ---
not arcs, since $\partial D$ is contained on $P$, and $\sigma$ acts as
a reflection in a small neighborhood of $P$.

Let $\gamma_1$ be a circle of intersection of $D$ and $\sigma(D)$
which is innermost on $\sigma(D)$.  Then $\gamma_1$ lies on $D$ and
$\sigma(D)$, hence $\gamma_2 = \sigma(\gamma_1)$ lies on $\sigma(D)$
and $D$, and $\gamma_2$ is innermost on $D$.  Surger: Replace $D$ by
replacing the disk bounded by $\gamma_1$ on $D$ with the disk bounded
by $\gamma_1$ on $\sigma(D)$, and push off $\sigma(D)$ slightly.  Call
this new disk $D'$.

We claim that the number of intersections $|D' \cap \sigma(D')|$ is
now less than $|D \cap \sigma(D)|$.  Outside a neighborhood of the
disk bounded by $\gamma_1$, $D'$ agrees with $D$.  Hence outside a
neighborhood of the disk bounded by $\gamma_2 = \sigma(\gamma_1)$,
$\sigma(D')$ agrees with $\sigma(D)$.

There are two cases to consider.  First, if $\gamma_2$ is outside the
disk bounded by $\gamma_1$ on $D$, then $\gamma_2$ and the disk it
bounds in $D$ are still contained in $D'$ (since $D'$ agrees with $D$
outside $\gamma_1$).  Similarly, $\gamma_1$ will be outside the disk
bounded by $\gamma_2$ on $\sigma(D)$, so $\gamma_1$ and the disk
bounded by $\gamma_1$ on $\sigma(D)$ will still remain on
$\sigma(D')$.  An example of this is illustrated in Figure
\ref{fig:P-incomp-1}.

\begin{figure}
\input{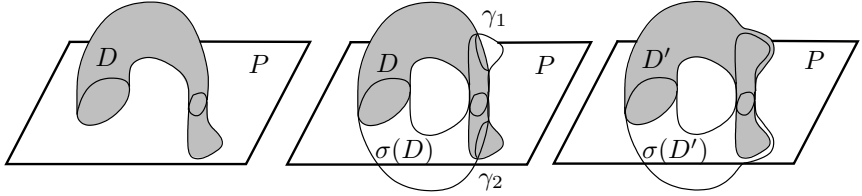}
\caption{Right: $D$ (shaded). Middle: $D$ and $\sigma(D)$ intersect at
	$\gamma_1$ and $\gamma_2$.  Left: Replace the interior of $\gamma_1$
	and push off to reduce the number of intersections by $2$.}
\label{fig:P-incomp-1}
\end{figure}

When we push off $\sigma(D)$ to form $D'$, we may do so equivariantly.
So $D'$ doesn't intersect $\sigma(D')$ in a neighborhood of the disk
bounded by $\gamma_1$ on $\sigma(D')$, and $\sigma(D')$ doesn't
intersect $D'$ in a neighborhood of the disk bounded by $\gamma_2$ on
$D'$.  Elsewhere, $D$ and $D'$ agree, so other intersections have not
changed.  Thus the number of intersections has decreased under this
operation.  

If instead $\gamma_2$ is inside the disk $E_1$ bounded by $\gamma_1$
on $D$, then when we surger, even before pushing off, $\gamma_1$ now
bounds a disk $E$ on $D'$, whose image under $\sigma$ is bounded by
$\gamma_2$ in $E_1$, hence is contained in $E_1$.  But now $E_1$ is
disjoint from $D'$.  Since $D$ and $D'$ agree elsewhere, $D'$ and
$\sigma(D')$ have fewer intersections.

Repeating this process a finite number of times, we obtain a
compressing disk $D$ such that $D$ and $\sigma(D)$ are disjoint.

Then $D \cup \sigma(D)$ is a sphere $S$ in $N$.  Since $N$ is
irreducible, $S$ bounds a ball $B$ in $N$ whose boundary $S$ is
invariant under $\sigma$.  Since $S$ meets $P$, so does $B$, and hence
$B$ must be preserved by the involution $\sigma$.

Now we have an orientation reversing involution of a ball $B$ which
fixes a circle $\partial D$ on the boundary of the ball and swaps the
disks on the boundary.  Double the ball across its boundary and extend
$\sigma$.  This gives an orientation reversing involution of $S^3$
with fixed point set a surface.  It follows from work of Smith in the
1930s that the fixed point set must be a $2$--sphere containing
$\partial D$, and $B$ intersected with this fixed point set must
therefore be a disk.  But the fixed point set of $\sigma$ is $P$, so
$P$ contains a disk with boundary $\partial D$.  This contradicts the
fact that $D$ was a compressing disk for $P$.
\end{proof}

We now use this result to show incompressibility of $D_i-L$.

\begin{lemma}
Let $L$ be a reduced generalized augmented link in $M$ such that $M-L$
is irreducible, with twisting disks $D_i$.  Then each $D_i-L$ is
incompressible in $M-L$.
\label{lemma:D-incompress}
\end{lemma}

\begin{proof}
Let $P$ be the surface of the reflection provided by Proposition
\ref{prop:reflect}, and suppose, by way of contradiction, that the
punctured disk $D_1-L$, say, is compressible.  Let $D$ be a
compressing disk in $M-L$.  Consider $D\cap P$.  Note $\D D \cap P$,
and therefore $D\cap P$, cannot be empty, else $\partial D$ bounds a
disk on $D_1-L$.  So $D\cap P$ consists of a non-empty collection of
arcs and closed curves on $D$.

We may assume $D\cap P$ does not contain any closed curves on $D$,
using the incompressibility of $P$, Lemma \ref{lemma:P-incompress}.
For if there is a closed curve component of $D\cap P$, it bounds a
disk on $P$, hence can be isotoped off.

Next we show that we can take $D$ to be invariant under $\sigma$.
Consider an outermost arc $\gamma$ of $D\cap P$.  This bounds a disk
$E$ whose boundary consists of the arc $\gamma$ and a portion of
$\partial D$.  Consider $F=E \cup \sigma(E)$.  This is a disk embedded
in $M-L$ with boundary on $D_1$, which is invariant under $\sigma$.
If $F$ is a compressing disk for $D_1 -L$, then replace $D$ by $F$,
and we have a compressing disk invariant under $\sigma$.  So suppose
$F$ is not a compressing disk for $D_1-L$.  Then $\D F$ bounds a disk
$E_1$ on $D_1$.  So $F \cup E_1$ is a 2--sphere, which must bound a
3--ball by irreducibility.  Isotope $D$ through this ball, pushing $E$
along $E_1$ through $P$, reducing the number of intersections of $P$
with $D$.  Since there are only finitely many intersections of $D\cap
P$, repeat this process only a finite number of times, and eventually
obtain a compressing disk for $D_1-L$ invariant under $\sigma$, or we
reduce to the case $D\cap P = \emptyset$, contradicting the fact that
$D$ is a compressing disk.

Now, $\partial D$ does not bound a disk on $D_1-L$.  Hence it must
bound a disk $E_1$ in $D_1$ punctured by $L$.  Replace $D_1$ by
replacing $E_1$ with $D$.  Then we have a new disk $D'$ which is
invariant under $\sigma$, which meets the $K_j$ fewer times.
But the generalized augmented link $M-L$ is reduced, which means, by
Definition \ref{def:reduced}, that the number of intersections $D_1
\cap (\cup K_j)$ is minimal over all disks in $M$ bounded by $C_1$.
Since $D'$ is another such disk, this is a contradiction.
\end{proof}

\subsection{Annuli}

Next, we show a series of results on annuli that are admitted in a
generalized augmented link.

\begin{lemma}[Lemma 2.5.3 of \cite{CGLS}]
  Let $M$ be an irreducible $3$--manifold with torus boundary
  components, not homeomorphic to $T^2\times I$.  Let $T_1$ and $T_2$
  be boundary components which are incompressible.  Suppose $A_1$ and
  $A_2$ are properly embedded annuli in $M$ with $\partial A_i =
  c_{i1} \cup c_{i2}$, with $c_{ij}$ on $T_j$.  Then $c_{1j}$ is
  isotopic to $c_{2j}$ on $T_j$, $j=1, 2$.
\label{lemma:distinct-annuli}
\end{lemma}

Lemma \ref{lemma:distinct-annuli} is actually not as general as
\cite[Lemma 2.5.3]{CGLS}.  Because the statement of that lemma is a
little different from Lemma \ref{lemma:distinct-annuli}, we reproduce
the proof here for convenience.

\begin{proof}
Assume $A_1$ and $A_2$ are in general position.  Let $\Delta_j$, $j=1,
2$ denote the number of intersections of $c_{1j}$ and $c_{2j}$.
Suppose one is nonzero.
  
There is an isotopy of the $A_i$ such that $\Delta_1 = \Delta_2$, and
any arc of intersection of $A_1 \cap A_2$ runs from one torus to the
other.  Otherwise, an arc of intersection runs from one torus back to
itself.  By an innermost arc argument, using the irreducibility of $M$
and the incompressibility of $T_i$, we can isotope $A_1$ and $A_2$ to
remove this arc of intersection. 

Now, we claim we can replace $A_2$, if necessary, so that $c_{1j}$ and
$c_{2j}$ intersect just once.  Suppose $c_{1j}$ and $c_{2j}$ intersect
at least twice.  Choose two arcs adjacent to each other on $A_2$, say
$a$ and $a'$.  That is, $a$ and $a'$ bound a disk $E_2$ on $A_2$ whose
interior is disjoint from $A_1$.  The arcs $a$ and $a'$ will also
bound a disk $E_1$ on $A_1$, whose interior is not necessarily
disjoint from $A_2$, but must be disjoint from $E_2$ by choice of
$E_2$.

The boundaries of $E_1$ are $a$, an arc on $c_{11}$ which we denote
$b_1$, $a'$, and an arc on $c_{12}$, which we denote $d_1$.
Similarly, the boundaries of $E_2$ are $a$, $b_2$ on $c_{21}$, $a'$,
and $d_2$ on $c_{22}$.  Then $E_1 \cup E_2$ gives an annulus $A$
embedded in $M$ with boundary components $b_1 \cup b_2$ on $T_1$ and
$d_1 \cup d_2$ on $T_2$.  Since ${\rm{int}}(b_2 \cap A_1) = \emptyset$
(because $E_2$ is disjoint from $A_1$), $A$ has slope $r_0$, say, on
$T_1$, where $\Delta(r_0, c_{i1}) = 1$.  Hence after isotopy, $A \cap
A_1$ will consist of a single arc (see \cite[Figure 2.3]{CGLS}).
Replace $A_2$ with $A$.
  
Now, $A_1 \cup A_2$ is homeomorphic to $X\times I$, where $X = X\times
\{0\}$ is the union of two simple loops on $T_1$ which intersect
transversely in a single point.  It has a regular neighborhood
homeomorphic to $N\times I$, where $N = N\times \{0\}$ is a regular
neighborhood of $X$ on $T_1$.  But now $\partial N$ bounds a disk
$D_1$ on $T_1$, $\partial N \times \{1\}$ bounds a disk $D_2$ on
$T_2$, so $D_1 \cup \partial N \times I \cup D_2$ is a 2--sphere,
which bounds a 3--ball $B$ in $M$ since $M$ is irreducible.  Then $M =
N\times I \cup B$ is homeomorphic to $T \times I$.
\end{proof}

\begin{corollary}
Suppose $M$ is an irreducible $3$--manifold with torus boundary
components, not homeomorphic to $T^2\times I$, which admits an
orientation reversing involution $\sigma$ which fixes a surface $P$
meeting incompressible components $T_1$ and $T_2$ of $\partial M$.
Suppose $A$ is an annulus embedded in $M$ with boundary components
lying on $T_1$ and $T_2$.  Then the slopes of $\partial A$ on $T_i$,
$i=1, 2$, are preserved by the involution $\sigma$.
\label{cor:fix-annulus}
\end{corollary}

\begin{proof}
  If $\sigma$ does not preserve one of the slopes $(\partial A)_i$
  on $T_i$, then $A$ and $\sigma(A)$ are two distinct annuli embedded
  in $M$ with non-isotopic boundary components.  This contradicts
  Lemma \ref{lemma:distinct-annuli}.
\end{proof}

We now apply these results to generalized augmented links.  First, we
need to rule out the case that a generalized augmented link might have
complement in $M$ homeomorphic to $T^2\times I$.


\begin{lemma}
Suppose $L$ is a generalized augmented link in $M$ and $M-L$ is
irreducible.  Then $M-L$ is not homeomorphic to $T^2\times I$.
\label{lemma:T2xI}
\end{lemma}

\begin{proof}
Suppose not.  Since $T^2\times I$ has just two boundary components,
and since we assume any generalized augmented link has at least one
crossing circle $C_1$, one boundary component of $M$ corresponds to
$C_1$ and the other to a knot strand.  Since the punctured $D_1$ is
incompressible by Lemma \ref{lemma:D-incompress} and $2$--sided, it
must be either horizontal or vertical in a Seifert fibering of
$T^2\times I$.  But then $D_1-L$ must be an annulus, contradicting the
fact that $m_1 \geq 2$.
\end{proof}

Using this fact, we can rule out annuli embedded in $M-L$.

\begin{lemma}
Let $L$ be a generalized augmented link in $M$ such that $M-L$ is
irreducible.  Then there is no annulus embedded in $M-L$ with one
boundary component on $\D\CN(C_i)$ parallel to $P \cap \D \CN(C_i)$,
and the other boundary component disjoint from $\D \CN(C_i)$.
\label{lemma:annulus-genaug}
\end{lemma}

\begin{proof}
Suppose $A$ is an annulus embedded in $M-L$ with boundary component
$(\D A)_i$ on $\D \CN(C_i)$ parallel to $P \cap \D \CN(C_i)$.

Consider the intersection of $D_i$ with $A$.  Since $(\partial A)_i$
is parallel to $P\cap \D \CN(C_i)$, and $\partial D_i$ meets $P$ twice
by Lemma \ref{lemma:meet-twice}, we may isotope $A$ so that $\partial
D_i \cap A$ consists of one or two points, depending on whether $P
\cap \partial \CN(C_i)$ consists of two or one components, as
illustrated on the left and right of Figure \ref{fig:half-twist},
respectively.

If $A\cap \partial D_i$ consists of one point, then $P\cap \partial
\CN(C_i)$ has two meridional components, and $(\partial A)_i$ is a
meridian on $\partial \CN(C_i)$.  But now consider $A\cap D_i$.  This
consists of arcs and curves of intersection.  Any arc has two
endpoints on $A \cap \partial D_i$, so there must be an even number of
points of intersection of $A\cap \partial D_i$.  However, we are
assuming there is just one such point.  This is a contradiction.

Thus $A\cap \partial D_i$ must consist of two points, and $A \cap D_i$
must have a single arc component running from $\partial D_i$ to
$\partial D_i$.  This arc bounds a disk $E$ in $D_i$, and a disk $E'$
in $A$, since the other boundary component of $A$ is disjoint from $\D
\CN(C_i)$.  We may assume the interiors of $E$ and $E'$ are disjoint
by an innermost curve argument, for intersections must be simple
closed curves in both, since $A\cap D_i$ has just one arc component,
and hence if $E$ and $E'$ are not disjoint we may isotope them off of
each other using irreducibility of $M-L$.  So $E \cup E'$ is a disk
with boundary on $\D \CN(C_i)$.  This must bound a disk on $\D
\CN(C_i)$.  By irreducibility of $M-L$, we may therefore isotope $A$
to have no intersections with $\D D_i$, contradicting the fact that
$(\D A)_i$ is parallel to $P \cap \D \CN(C_i)$.
\end{proof}

\begin{lemma}
If $L$ is a reduced generalized augmented link in $M$ such that $M-L$
is irreducible, then there is no annulus embedded in $M-L$ with
boundary components on $\partial \CN(C_i)$, $\partial \CN(C_j)$, for
$i\neq j$.
\label{lemma:annulus}
\end{lemma}

\begin{proof}
By Definition \ref{def:reduced}(1), each $\D \CN(C_i)$ is
incompressible.  Thus by Corollary \ref{cor:fix-annulus}, any embedded
annulus $A$ must have boundary components fixed by $\sigma$.  At most
two slopes on $\D \CN(C_i)$ are fixed by $\sigma$.  These are the
slopes of $P \cap \D \CN(C_i)$ and of $D_i \cap \D \CN(C_i) = \D D_i$.

By Lemma \ref{lemma:annulus-genaug}, no boundary component of $A$ is
parallel to $P \cap \D \CN(C_i)$ or to $P \cap \D \CN(C_j)$.  By
Definition \ref{def:reduced}(2), we cannot have $(\D A)_i$ and $(\D
A)_j$ parallel to $\D D_i$ and $\D D_j$.  Thus no such annulus exists.
\end{proof}


\subsection{Seifert fibered augmented links}

We may now classify all Seifert fibered reduced generalized augmented
links.

\begin{theorem}
If $M-L$ is irreducible and Seifert fibered, where $L$ is a reduced
generalized augmented link in $M$, then $L$ has just one crossing
circle component $C_1$, $M-C_1$ is a solid torus, $P$ is an embedded
annulus or M\"obius band in $M-C_1$, and the knot strands are parallel
to the boundary of $P$.  In particular, if $P$ is an annulus, there
are at least two knot strand components.
\label{thm:no-sf}
\end{theorem}

\begin{figure}
\begin{center}
\input{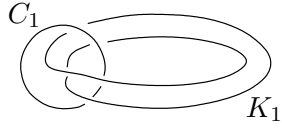}
\end{center}
\caption{The only Seifert fibered reflective augmented link with one
	knot strand component.}
\label{fig:mobius}
\end{figure}

\begin{proof}
Suppose $M-L$ is Seifert fibered. First, by Lemma \ref{lemma:annulus},
there can be no annuli between link components $C_i$ and $C_j$.  This
implies that there cannot be more than one link component $C_1$.

Now, $D_1-L$ is incompressible by Lemma \ref{lemma:D-incompress}.
Since $D_1-L$ is $2$--sided, it is horizontal or vertical in $M-L$.
If vertical, it must be an annulus, contradicting the fact that $m_1
\geq 2$.  So $D_1-L$ is horizontal.  Then the meridians of the knot
strands (i.e. the curves which bound disks in $M$) cannot be Seifert
fibers, so the Seifert fibering of $M-L$ extends to $M-C_1$.  The base
orbifold of $M-C_1$ is branch covered by the horizontal surface $D_1$,
hence it is a disk with one singular point. Thus $M-C_1$ must be a
solid torus.

Since $P$ is incompressible by Lemma \ref{lemma:P-incompress}, if it
is orientable, then it is an annulus.  Because the knot strands are
embedded in $P$ and nontrivial, they must be parallel to the core of
the solid torus $M-C_1$.  In particular, if there is just one knot
strand component, then $M-L$ is homeomorphic to $T^2 \times I$,
contradicting Lemma \ref{lemma:T2xI}.  So in this case there are at
least $2$ knot strand components.

If $P$ is non-orientable, then by work of Frohman
\cite{frohman:pseudo-vert} and Rannard \cite{rannard:pseudo-vert}, $P$
is \emph{pseudo-vertical} in a solid torus, meaning, in this case, it
is a punctured non-orientable surface in the solid torus $M-C_1$.
These were classified by Tsau \cite{tsau:incompr-solid-tori}, and have
boundary of the form $P \cap \partial \CN(C_1) = \alpha =q \mu + (2k)
\lambda$, where $\mu$ is a meridian of the solid torus, $\lambda$ is a
longitude, $k\geq 1$, and $q$ is an odd integer.

In our case, we know which boundary slopes $\alpha$ can occur, because
of the existence of the involution $\sigma$.  In particular, since
$\partial D_1$ intersects $P$ exactly twice, by Lemma
\ref{lemma:meet-twice}, $k=1$.  Then by untwisting, we may assume
$q=1$, and so $\alpha$ is the slope $\mu + 2\lambda$, and $P$ is a
M\"obius band in the solid torus $M-C_1$.

The knot strands are embedded in the M\"obius band $P$ and nontrivial.
Thus they must be parallel to $\D P$.  If there is just one knot
strand, the link $L$ is as in Figure \ref{fig:mobius}.  More
precisely, it has complement homeomorphic to the complement of the
link of Figure \ref{fig:mobius} in $S^3$.
\end{proof}

We obtain the following immediate consequence of this result.

\begin{corollary}
Let $K$ be a knot in $S^3$ which has a diagram whose augmentation is a
Seifert fibered reduced generalized augmented link.  Then $K$ is a
$(2,q)$ torus knot.
\label{cor:sf}
\end{corollary}

\begin{proof}
Let $L$ denote the augmentation of $K$.  First we show $S^3-L$ is
irreducible.  The manifold $S^3-L$ is homeomorphic to $(S^3-K)-(\cup\,
C_i)$, where $C_i$ are crossing circles encircling generalized twist
regions.  Let $S$ be a sphere in $S^3-L$, and assume it does not bound
a ball in $S^3-L$.  Then the image of $S$ under the homeomorphism is a
sphere in $(S^3-K)-(\cup\, C_i)$ which does not bound a ball.  Since
$S^3-K$ is irreducible, $S$ must bound a ball $B$ in $S^3-K$.  Hence
some $C_i$ lies in $B$.  But $C_i$ is unknotted, hence bounds a disk
in $B$.  This contradicts property (1) of the definition of reduced,
Definition \ref{def:reduced}.

By Theorem \ref{thm:no-sf}, the diagram $D$ can have only one
generalized twist region.  The augmentation is the link shown in
Figure \ref{fig:mobius}, since there is just one component $K$.  Thus
when we twist to obtain $S^3-K$, we remove $C_1$ from the diagram and
add an even number of crossings at the twist region it bounds.  This
is a $(2,q)$ torus knot.
\end{proof}

\section{Essential tori and augmented links}
\label{sec:satellite}
In this section, we consider reduced generalized augmented links $L$
such that $M-L$ is toroidal.  We show that the torus decomposition of
$M-L$ satisfies some nice properties.

Recall that by work of Jaco and Shalen \cite{jaco-shalen} and
Johannson \cite{johannson}, every irreducible 3--manifold $N$ with
(possibly empty) torus boundary contains a pairwise disjoint
collection of embedded essential tori $\TT$, unique up to isotopy,
such that the closure of a component of $N-\TT$ is either atoroidal or
Seifert fibered.  If $N$ admits an orientation reversing involution
$\sigma$, then by the equivariant torus theorem, first proved by
Holzmann \cite{holzmann}, each incompressible torus in $N$ is isotopic
to one which is preserved by $\sigma$ or taken off itself.  Then
$\sigma$ applied to $\TT$ gives a new torus decomposition of $N$,
which by uniqueness must agree with $\TT$.  Thus the closure of each
component of $N-\TT$ is either fixed by $\sigma$, or taken off itself.
This is the equivariant torus decomposition of Bonahon and Siebenmann
\cite{bon-sieb}.  We refer to this equivariant torus decomposition as
the JSJ decomposition.

\begin{theorem}
Let $L$ be a reduced generalized augmented link in $M$, such that
$M-L$ is irreducible.  Then there exists a collection $\mathcal{T}$ of
tori, incompressible in $M-L$, such that if $U$ is a component of
$M-\mathcal{T}$, $U$ satisfies one of the following:
\begin{itemize}
	\item $U$ does not meet any crossing circles of $L$.
	\item $U$ meets exactly one crossing circle $C_i$, $U$ is
	homeomorphic to $T^2\times I$, and $C_i$ is isotopic to a simple
	closed curve on $\D U$.
	\item $U-(L\cap U)$ is the complement of a reduced generalized
	augmented link.
\end{itemize}
\label{thm:tori}
\end{theorem}

The collection of tori $\mathcal{T}$ may be larger than the minimal
collection of the JSJ decomposition.  However, we find $\mathcal{T}$
by adding incompressible tori to the JSJ decomposition.

\begin{lemma}
Let $L$ be a reduced generalized augmented link in $M$ such that $M-L$ is
irreducible.  Let $\TT$ denote the tori of the JSJ decomposition for
$M-L$.  Suppose there is a component $U$ of $M-\TT$ that contains a
crossing circle $C_j$ and an embedded annulus $A$ with one boundary
component on some $T_0 \subset \D U$ and one on $C_j$.  Let $T_j$ be
the torus obtained by taking the boundary of a small regular
neighborhood of the union of $T_0$, $A$, and $C_j$.  Then $T_j$ is
incompressible in $M-L$.
\label{lemma:finding-Tj}
\end{lemma}

The collection $\mathcal{T}$ of Theorem \ref{thm:tori} will be
obtained by adding to $\TT$ any $T_j$ of Lemma \ref{lemma:finding-Tj}
which are not isotopic to tori in $\TT$.

\begin{figure}
\input{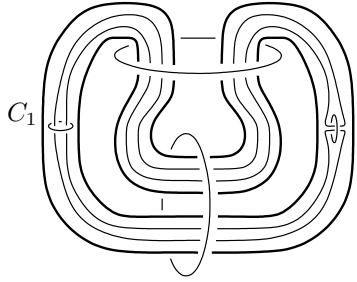}
\caption{$C_1$ and the torus denoted by the thick line bound an
	annulus.}
\label{fig:toroidal-link}
\end{figure}

Before proving the lemma, we illustrate by example a situation in
which the lemma will apply.  Consider the link in Figure
\ref{fig:toroidal-link}.  The heavy line in that figure shows the
location of an incompressible torus $T_0$.  When we cut along $T_0$,
we obtain a hyperbolic generalized augmented link on the outside,
homeomorphic to the complement of the Borromean rings.  On the inside,
$C_1$ and $T_0$ bound an annulus.  Take the boundary of a regular
neighborhood of the union of this annulus with $T_0$ and $C_1$, and we
obtain an incompressible torus $T_1$ as in Lemma
\ref{lemma:finding-Tj}.  Cutting along $T_1$, we split the inside into
two components, one homeomorphic to $(T^2\times I)-C_1$, and the other
(in this example) another copy of the Borromean rings complement.

Recall that we are interested in links obtained by twisting
generalized augmented links.  When we do twisting along $C_1$ in the
above example, the component $(T^2\times I)-C_1$ becomes the manifold
$T^2 \times I$.  Thus the two incompressible tori $T_0$ and $T_1$
become isotopic to each other after twisting.  By the results in
\cite{purcell:aug}, sufficiently high twisting along the remaining
crossing circles in each gives a manifold which remains hyperbolic.
Thus the torus decomposition of the twisted link contains
just one of $T_0$ and $T_1$.  We will see this in Theorem
\ref{thm:dehnfill}.

\begin{proof}[Proof of Lemma \ref{lemma:finding-Tj}]
Suppose by way of contradiction that $T_j$ is compressible.  A
compressing disk can be isotoped to lie in $U-(L\cap U)$, else we
obtain a compressing disk for the torus $T_0$ in $\TT$, which is
impossible.

Surger along a compressing disk for $T_j$ in $U-(L \cap U)$ to obtain
a sphere $S$ embedded in $U-(L \cap U)$.  Note that $U-(L \cap U)$ is
irreducible, since $M-L$ is irreducible and $T_0$ is incompressible.
Thus $S$ bounds a ball $B$ in $U-(L\cap U)$.  Now, $B$ cannot be on
the side of $S$ containing $A$ since this side contains boundary
components $T_0$ and $C_j$ of $U-(L \cap U)$.  Thus $T_j$ bounds a
solid torus $V$ in $U-(L \cap U)$.  Recall that $T_j$ was formed by
taking the boundary of a regular neighborhood of the union of $T_0$,
an annulus $A$, and $C_j$ in $U$.  Since $T_j$ bounds a solid torus
$V$ in $U-(L\cap U)$, it must be the case that $U-(L \cap U)$ is
homeomorphic to $V$ union a regular neighborhood of the annulus $A$.
This has just two boundary components: $T_0$ and $C_j$.  We claim this
is impossible.

Since $U$ contains $C_j$, the surface $P$ of Proposition
\ref{prop:reflect} meets $U$, and the torus boundary components of $U$
are preserved by $\sigma$.  Hence $U-(\cup\, C_i)$ is preserved by
$\sigma$, where the union is over $C_i$ in $U$.  The annulus $A$ has
one boundary component, $A_1$, say, on $\D \CN(C_j)$, taken by
$\sigma$ to $-A_1$.  Hence it is isotopic to an annulus which meets
$P$ in two arcs and is preserved under $\sigma$.  Since tori $T_0$ and
$\partial \CN(C_j)$ are also preserved under $\sigma$, the solid torus
$V$ is preserved under $\sigma$.  Then $V \cap P$ must be an annulus,
a M\"obius band, or two meridional disks in $V$.

We form $U-(L\cap U)$ by attaching a thickened annulus to $V$.  This
thickened annulus is attached along some slope $\mu$ on $\partial V$.
Since $A$ is taken to itself with reversed orientation by $\sigma$,
the slope $\mu$ must be taken to $-\mu$ by $\sigma$.  There are very
few possibilities for $\mu$.

In case $V \cap P$ is an annulus or M\"obius band, $\mu$ must bound a
disk in $V$.  Attaching an annulus to $V$ along two meridians gives a
manifold with compressible boundary, but neither $C_j$ nor $T_0$ is
compressible.

Thus $V \cap P$ consists of two meridional disks, and $\mu$ must be
some longitude of $\partial V$.  When we attach a thickened annulus to
longitude slopes, the resulting manifold is homeomorphic to $T^2
\times I$.  Then $T_0$ is parallel to $C_j$, contradicting the fact
that $T_0$ is essential in $M-L$.
\end{proof}

\begin{lemma}
Let $C_i$ and $C_j$, $i\neq j$, be crossing circles that satisfy the
hypotheses of Lemma \ref{lemma:finding-Tj}.  Then the tori $T_i$ and
$T_j$ are not isotopic.
\label{lemma:ti-tj}
\end{lemma}

\begin{proof}
$T_i$ is formed by taking the boundary of a regular neighborhood of
the union of an incompressible torus, an annulus with boundary on the
torus and on $C_i$, and $C_i$.  Thus there is an annulus $A_i$ with
boundary on $T_i$ and on $C_i$.  Similarly for $T_j$ there is an
annulus $A_j$ with boundary on $T_j$ and on $C_j$.  If $T_i$ and $T_j$
are isotopic, then the annulus $A_j$ may be isotoped to have boundary
on $C_j$ and on $T_i$.  Then we form a new annulus between $C_i$ and
$C_j$ by taking the union of $A_i$, $A_j$, and an annulus running
along $T_i$ between these two.  This contradicts the fact that $L$ is
reduced.
\end{proof}

Form the collection $\mathcal{T}$ of Theorem \ref{thm:tori} by adding
to $\TT$ all the tori $T_j$ of Lemma \ref{lemma:finding-Tj} which are
not isotopic to a torus in $\TT$.  By Lemma \ref{lemma:ti-tj}, these
are not isotopic to each other.  The closure of each component of
$(M-L) - \mathcal{T}$ is still either atoroidal or Seifert fibered,
but $\mathcal{T}$ may no longer be the minimal such collection.  By
construction, components of $M - \mathcal{T}$ either do not contain
crossing circles of $L$; contain a single crossing circle $C_j$
sandwiched between incompressible tori $T_0$ and $T_j$ of Lemma
\ref{lemma:finding-Tj}, in which case the second possibility of
Theorem \ref{thm:tori} holds; or $U$ contains crossing circles of $L$,
but none of these bound annuli with boundary on $\D U$.  For the proof
of Theorem \ref{thm:tori}, we need to examine these components, and
show that any such $U-(L\cap U)$ is homeomorphic to the complement of
a reduced generalized augmented link.

Generalized augmented links are defined to lie in a 3--manifold
admitting a reflection through a surface $S$.  First, we define a
manifold $N$ which will play the role of thus underlying 3--manifold.

Let $U$ be a component of $M-\mathcal{T}$ which contains at least one
crossing circle $C_i$, but does not contain an embedded annulus with
boundary on $\D U$ and on $C_i$, for any $C_i \subset U$.  Consider $D_i
\cap U$.  In the manifold $M$, $D_i$ may intersect essential tori of
$\mathcal{T}$.  Then $D_i$ will meet boundary components of $U$.  Let
$\ell_{i_1}, \dots, \ell_{i_k}$ be the boundary components of that
component of $D_i\cap U$ which meets $C_i$.  Replace $U$ by the
manifold $N_i$ obtained by Dehn filling $V$ along the slopes
$\ell_{i_j}$, $j=1, \dots, k$.  Let $\hat{K}_{i_1}, \dots,
\hat{K}_{i_k}$ denote the solid tori attached in the Dehn filling.

Do this Dehn filling for each $C_i$ contained in $U$. We obtain a new
manifold $N$.  Note $N$ is well--defined because the $D_i$ are
disjoint; thus if $T_k$ is met by $D_i$ and $D_j$, then $D_i \cap T_k$
and $D_j \cap T_k$ must give the same slope, so the Dehn fillings
along that slope are the same.  Similarly, $D_i$ might meet $T_k$
several times, but again along the same slope.


Now, let $L_U$ be the link in $N$ consisting of components $C_i\cap U$
and $K_j\cap U$, as well as the cores of each distinct solid torus in
the set $\{\hat{K}_{k}\}$.  We will abuse notation slightly and
continue to refer to these cores of solid tori by $\hat{K}_k$.  The
following is immediate.

\begin{lemma}
For $N$ the manifold and $L_U$ the link in $N$ constructed as above,
$N-L_U$ is homeomorphic to $U-(L\cap U)$. \qed
\label{lemma:Nhomeo}
\end{lemma}

We claim that $L_U$ is a reduced generalized augmented link in $N$,
with components $C_i\cap U$ taking the role of the crossing circles,
and components $K_j\cap U$ and $\hat{K}_k$ taking the role of the knot
strands.

By assumption, there is at least one $C_i$, bounding $D_i$, meeting at
least one $K_i$ or $\hat{K}_i$.  So the link $L_U$ contains at least
the minimal number of necessary link components to be a generalized
augmented link.

\begin{lemma}
The manifold $N$ admits an involution through a surface $S\subset N$,
a link $K$ is contained in a neighborhood of $S$, has diagram $D$ and
maximal twist region selection such that when we encircle generalized
twist regions of $D$ by crossing circles and untwist, the result is a
link isotopic to $L_U$.  That is, $L_U$ is a generalized augmented
link, given by augmenting a diagram of $K$ in $N$.
\label{lemma:L_U-aug}
\end{lemma}

\begin{proof}
The involution $\sigma$ of Proposition \ref{prop:reflect} preserves
$U-(\cup\, C_i)$, where the union is over crossing circles $C_i$ in
$U$.  It has fixed point set $U\cap P$, and components $K_j$ in $L_U$
are embedded in $P$.  We show that the involution $\sigma$ extends to
the solid tori $\CN(\hat{K}_k)$, and that the cores $\hat{K}_k$ are
embedded in the surface $P$.

Recall that $\CN(\hat{K}_k)$ has boundary which is an incompressible
torus $T_k$ in $M-L$, and $T_k$ is preserved by $\sigma$.  Moreover,
some $D_j$ meets $T_k$ in a meridian of $\CN(\hat{K}_k)$.  The slope
$\partial D_j \cap T_k$ cannot bound a disk in $M-L$ by
incompressibility of $T_k$.  Since it does bound a disk in $D_j$, this
disk must be punctured by some component $K_i$ in $M-L$.  Therefore
the slope $\D D_j \cap T_k$ must be taken by $\sigma$ to $-\D D_j \cap
T_k$.  This means a meridian of the solid torus $\CN(\hat{K}_k)$ is
inverted by the involution $\sigma$.  Since the boundary is preserved
by $\sigma$, it follows that the involution extends to give an
involution of the solid torus $\CN(\hat{K}_k)$.  Also, $P \cap
\CN(\hat{K}_k)$ must be a longitude of $\partial \CN(\hat{K}_k)$, in
the sense that it intersects the meridian exactly once.  Therefore the
core, $\hat{K}_k$, is embedded in $P$.

Now note that $N$ and $L_U$ satisfy the conclusion of Proposition
\ref{prop:reflect}.  Since $\sigma$ acts on $N$ as an extension of the
involution of Proposition \ref{prop:reflect} acting on $M-(\cup\,
C_i)$, we know that $P$ meets the $C_i$ in $N$ in the same way it
meets the $C_i$ in $M$.  Namely, $P$ either meets $\CN(C_i)$ in two
meridian components, as on the left in Figure \ref{fig:half-twist}, or
in a single component as on the right in Figure \ref{fig:half-twist}.
Form the surface $S$ by taking $\hat{S}$ to be $P$ outside of a
neighborhood of those crossing disks that meet half--twists.  Inside a
neighborhood of a half--twist, $P$ will appear as on the right of
Figure \ref{fig:half-twist}.  The surface $\hat{S}$, however, should
run straight through the crossing circle, meeting $\CN(C_i)$ in two
meridians, as on the left of Figure \ref{fig:half-twist}.  Note there
is a reflection $\tau$ in $\hat{S}$ which still preserves $N-(\cup\,
C_i)$, although it reverses crossings at half--twists.  Moreover, the
reflection $\tau$ preserves a meridian of each crossing circle, hence
extends to a reflection $\tau$ of $N$ through the surface $S$, where
$S$ is obtained from $\hat{S}$ by capping off boundary components on
the $\CN(C_i)$ by disks.

Notice that the strands $K_i$ and $\hat{K}_k$ lie in a neighborhood of
$S$, as do the crossing circles $C_i$.  Hence when we twist along
crossing circles, we form a link $K$ which still lies in a
neighborhood of $S$.  Moreover, we may project $K$ to $S$ such that
the twists obtained by twisting along the $C_i$ form distinct
generalized twist regions, and so that $L_U$ is an augmentation of a
diagram of $K$.
\end{proof}

We now need to show that $L_U$ is reduced in $N$.  To do so, we find
disjoint embedded disks bounded by the crossing circles $C_j$ in $N$,
and show these are minimal in the sense of part (1) of Definition
\ref{def:reduced}.  

\begin{lemma}
Let $C_{i_1}, \dots, C_{i_k}$ be the crossing circles of $L_U$.  Each
bounds a disk $D'_{i_j}$ in $N$ such that the collection
$\{D'_{i_j}\}$ is embedded in $N$, the involution $\sigma$ restricts
to an orientation reversing involution of each $D'_{i_j}$, and the
disks meet $K_i$ and $\hat{K}_i$ in $m_i$ points, where $m_i \geq 2$
and $m_i$ is minimal over all disks in $M$ bounded by $C_i$.
\label{lemma:exist-disks}
\end{lemma}

\begin{proof}
Let $C_{i_1}, \dots, C_{i_k}$ be the crossing circles in $U$.  Each
$C_{i_j}$ bounds a disk $D_{i_j}$ in $M$.  By construction of $N$, the
collection $\{D_{i_1}, \dots, D_{i_k}\}$ extends to an embedded
collection of disks $\{D'_{i_1}, \dots, D'_{i_k}\}$ in $N$, such that
$\sigma$ restricts to an orientation reversing involution of each
$D'_{i_j}$.

If the $m_i$ is not minimal, then we can find a collection of
disjoint, embedded punctured disks which meet $\partial U$ fewer
times.  Replace the collection $\{D'_{i_1}, \dots, D'_{i_k}\}$ with
this new collection.  Each must meet $K_i$ and $\hat{K}_i$ in at least
$2$ points, for otherwise we would have an annulus between some $C_j$
and a component $K_i$ or $\hat{K_i}$.  The first cannot happen by
definition of a reduced generalized augmented link in $M$.  The second
cannot happen by assumption: $U$ is assumed to be a component of
$M-\CT$ which does not contain an embedded annulus with boundary
components on $C_i$ and on $\D U$, and $\D \CN \hat{K_i}$ is a
component of $\D U$.
\end{proof}

We are now ready to complete the proof of Theorem \ref{thm:tori}.

\begin{proof}[Proof of Theorem \ref{thm:tori}]
Let $\CT$ be the collection of tori described after the proof of Lemma
\ref{lemma:ti-tj}.  Let $U$ be a component of $M-\CT$.  If $U$ does
not contain any crossing circle $C_i$, then we are done.  If $U$
contains a crossing circle $C_i$ and an embedded annulus with boundary
$C_i$ and boundary on $\D U$, then by construction of $\CT$, $U$ is
homeomorphic to $T^2\times I$, $C_i$ is the only crossing circle
contained in $U$, and the curve $C_i$ is boundary parallel in $U$.
This is the second case of the theorem.

So assume $U$ contains a crossing circle, but does not contain any
embedded annuli with boundary on $\D U$ and on $C_i$.  Then by Lemma
\ref{lemma:Nhomeo}, $U-(L\cap U)$ is homeomorphic to the manifold
$N-L_U$, which, by Lemma \ref{lemma:L_U-aug} is a generalized
augmented link complement.  By Lemma \ref{lemma:exist-disks}, $L_U$
satisfies condition (1) of the definition of reduced, Definition
\ref{def:reduced}.  It satisfies condition (2) as well, since any
annulus embedded in $N-L_U$ with boundary components on $C_i$ and
$C_j$ is embedded in $M-L$ with boundary components on $C_i$ and
$C_j$.  Since $L$ is reduced, no such annulus exists.  The link $L_U$
satisfies condition (3) of Definition \ref{def:reduced} by assumption,
given the definition of $\CT$.
\end{proof}

\begin{theorem}
Let $L$ be a reduced generalized augmented link in $M$ with $M-L$
irreducible, and let $\CT$ be the tori of Theorem \ref{thm:tori}.  Let
$K$ be the link formed by twisting along all the crossing circles of
$L$, subject to the restriction that if $C_i$ is contained in a
component of $M-\CT$ which is not homeomorphic to $T^2\times I$, then
at least $6$ half--twists are inserted when we twist along $C_i$.
Then there is a torus decomposition of $M-K$ for which components of
the decomposition
\begin{itemize}
	\item are either atoroidal or Seifert fibered,
	\item are in one-to-one correspondence with the components of $M-\CT$
which are not homeomorphic to $T^2\times I$,
  \item and have the same geometric type (hyperbolic or Seifert
fibered) as the corresponding component of $M-\CT$.
\end{itemize}
\label{thm:dehnfill}
\end{theorem}


\begin{remark}
The decomposition of Theorem \ref{thm:dehnfill} may not be the JSJ
decomposition.  In particular, there may be two Seifert fibered
components which our decomposition separates, but which are considered
as one in the minimal JSJ decomposition.  We will see in Section
\ref{sec:knot} that this does not happen when $L$ comes from the
augmentation of a knot in $S^3$, but it could happen more generally.
\end{remark}

\begin{proof}
Let $U$ be a component of $M-\CT$.  If $U$ contains no crossing
circles, then twisting does not affect $U$ and so the result holds.

Similarly, if $U$ is homeomorphic to $T^2\times I$, and contains just
a single boundary parallel $C_i$, then twisting yields a manifold
homeomorphic to $T^2\times I$, which will not be a component of a
torus decomposition.

If $U$ contains a crossing circle, but is not $T^2\times I$, then
$U-(L\cap U)$ is either hyperbolic or Seifert fibered.

If it is hyperbolic, by \cite[Proposition 3.5]{purcell:aug}, the slope
of the twisting on a horoball neighborhood of a crossing circle $C_i$
has length at least $\sqrt{(1/4)+c_i^2}$, where $c_i$ is the number of
half--twists inserted.  Since $c_i\geq 6$, this length is greater than
$6$, hence by the 6--Theorem \cite{agol:bounds, lackenby:word}, the
result of Dehn filling is hyperbolic.

If $U-(L\cap U)$ is Seifert fibered, Theorem \ref{thm:no-sf}
implies $U-(L\cap U)$ is homeomorphic to the complement of parallel
strands embedded in an annulus or M\"obius band in a solid torus.
After twisting, we obtain a $(k, k)$ torus link, or a $(2k, pk)$ torus
link, where $k$ is the number of knot strands of $L_U$.  These torus
links are still Seifert fibered.
\end{proof}


\section{Reducing knot diagrams}
\label{sec:diagram}
We wish to apply the previous results to as many knots and links as
possible.

In this section, we prove that all knots in $S^3$ admit a diagram such
that the augmentation is reduced, as in Definition \ref{def:reduced}.
We say the diagram $D$ of a link $K$ is \emph{twist reduced} if there
exists a maximal twist region selection such that the corresponding
augmentation of $K$ gives a reduced generalized augmented link.

\begin{theorem}
Let $K$ be a knot in $S^3$ with diagram $D$ and a maximal twist region
selection.  Then there exists a twist reduced diagram $D'$ for $K$.
\label{thm:reduced-knot}
\end{theorem}

We will find the diagram of Theorem \ref{thm:reduced-knot} by forming
the augmentation of the given diagram, and then removing unnecessary
crossing circles and extracting unnecessary knot strands from crossing
disks.  When we do twisting on remaining crossing circles, projecting
twists to the projection plane in $S^3$ in the usual way, we will
obtain the desired diagram of the theorem.

\begin{define}
A \emph{standard diagram} of a generalized augmented link is a diagram
such that all knot strands lie on the projection plane except at
half--twists, which are contained in a neighborhood of the
corresponding crossing circle.  Crossing circles are perpendicular to
the projection plane, and crossing disks project to straight lines
running directly under the crossing circles of the diagram.  For
example, the portions of the diagrams in Figures
\ref{fig:cross-cir}(b), \ref{fig:nugatory}, and \ref{fig:redundant}
are standard.
\label{def:standard-diagram}
\end{define}

The next few lemmas ensure part (1) of Definition \ref{def:reduced}
will hold.

\begin{lemma}
Let $L$ be a generalized augmented link in $S^3$ with standard
diagram.  Suppose there exists a disk $E$ embedded in $S^3$, with
boundary some crossing circle $C_i$, disjoint from the other crossing
circles, and suppose that $E$ meets the knot strands fewer times than
does $D_i$.  Then there exists such a disk $F$ such that in addition,
$F \cap (\cup\, {\rm int}(D_j)) = \emptyset$.
\label{lemma:meetnoDi}
\end{lemma}

\begin{proof}
Suppose $E$ meets the interior of some $D_j$.  We may assume the
intersection is transverse and consists of simple closed curve
components.  There is some innermost disk $\hat{E}$ on $E$ whose
boundary is a curve $\gamma$ on $D_j$.  Consider the disk $G$
constructed by taking the disk $D_j$ outside $\gamma$, replacing $D_j$
inside $\gamma$ by $\hat{E}$.  Push $G$ off $D_j$ slightly, so $G$ and
$D_j$ do not intersect.

Suppose first that $G$ meets the knot strands fewer times than does
$D_j$.  Then replace $E$ by $G$, replacing $C_i$ by $C_j$.  This disk
$G$ has fewer intersections with $\cup\, {\rm int}(D_k)$ than does
$E$.

Suppose instead $G$ meets the knot strands the same number of times as
does $D_j$.  Replace $E$ by replacing $\hat{E} \subset E$ with the
portion of $D_j$ bounded by $\gamma$, and push off $D_j$.  We have
decreased the number of intersections of $E$ with $D_j$ without
increasing the number of intersections of $E$ with knot strands.

In either case, we have a new disk which meets the interiors of the
$D_k$ fewer times.  Repeat a finite number of times, and we obtain a
disk $F$ as in the statement of the lemma.
\end{proof}

\begin{lemma}
Let $L$ be a generalized augmented link in $S^3$ with standard
diagram.  Suppose there is a disk $E$ embedded in $S^3$ such that $E$
has boundary some crossing circle $C_i$, is disjoint from the other
crossing circles, and $E \cap (\cup\, {\rm int}(D_k))$ is empty.  Then
there exists such a disk $F$ such that in addition, $F$ intersects the
projection plane in a single arc $\gamma_1$, the intersections of $F$
with knot strands all lie on $\gamma_1$, and the number of intersections
of $F$ with knot strands is at most the number of intersections of $E$
with the knot strands.
\label{lemma:SintersectP}
\end{lemma}

\begin{proof}
First, we show we may assume that $E$ does not meet the diagram of $L$
in the neighborhood of any $D_j$ containing a half--twist.  That is,
we show we may assume $E$ does not meet any half--twists of the
diagram.  Assume $E$ does run through a half--twist corresponding to
$D_j$.  The half--twist is contained in a neighborhood of $D_j$, which
is homeomorphic to $D_j \times [-1,1]$.  Without loss of generality,
assume the half--twist is contained in $D_j \times (-1, 0)$, with
$D_j$ lying at $D_j \times \{0\}$.  Since $E$ does not meet $D_j$, it
must intersect $\D (D_j \times (-1,0))$ in the surfaces $D_j\times
\{-1\}$ or $\D D_j \times (-1, 0)$.  If any intersections of $E$ with
$\D (D_j \times (-1,0))$ bound disks in $\D (D_j \times (-1,0))-L$,
then we may replace $E$ by replacing corresponding disks in $E$ with
those in $\D (D_j \times (-1,0))-L$, and pushing into $D_j \times
(-1,0)$ slightly.  After this replacement, we may assume that any
component of $E \cap \D(D_j \times (-1,0))$ lies on $D_j \times
\{-1\}$ and bounds punctures of $D_j - L$.  Since $(D_j \times
(-1,0))-L$ is homeomorphic to $(D_j-L)\times (-1,0)$, we may replace
$E$ by replacing a disk bounded by $E\cap (D_j\times \{-1\})$ in $E$
by the corresponding disk in $D_j$, and pushing out of $D_j \times
(-1,0)$ slightly.  This will meet the knot strands of $L$ at most as
many times as does $E$, and will not intersect the neighborhood of
$D_j$ containing the half--twist.

So assume $E$ does not meet the diagram of $L$ in any half--twists.
Since outside of half--twists, all knot strands lie on the projection
plane, any intersections of $E$ with the knot strands must lie on the
projection plane.  It remains to show we can assume there is just one
component of intersection of $E$ with the projection plane.

Let $\gamma_1$ be the arc of intersection of $E$ with the projection
plane whose endpoints lie on $C_i$.  Consider $S^3$ cut along the
projection plane.  Remove neighborhoods of all the $D_j$ for $j\neq
i$, including all half--twists, as well as a neighborhood of a
half--twist at $D_i$, if applicable, which does not contain $C_i$.
This gives two balls, with $\gamma_1$ embedded in the boundary of each
ball, and one half of $C_i$ embedded as an arc in each ball, with
endpoints meeting those of $\gamma_1$.  Note no components of $L$
intersect the interior of either ball, aside from $C_i$.  Thus in each
ball we may find an embedded disk with boundary running along $C_i$
and along the arc $\gamma_1$ which only meets $L$ in its boundary.
These disks glue to give a new disk $F$ which meets the knot strands
exactly where $\gamma_1$ meets the knot strands, hence meets the knot
strands at most as many times as does $E$.  The disk $F$ satisfies the
conclusions of the lemma.
\end{proof}

\begin{lemma}
Let $L$ be a generalized augmented link in $S^3$ with standard
diagram.  Suppose there is a disk $E$ in $S^3$ with boundary some
crossing circle $C_i$, disjoint from the other crossing circles, such
that $E$ meets the knot strands fewer times than does $D_i$, and $E
\cap (\cup\,{\rm D_j})$ is empty.  Then we may isotope $L$ to a
generalized augmented link with standard diagram with the same
crossing circles, and the same number of strands running through each
crossing circle, except that $D_i$ meets the knot strands fewer times.
\label{lemma:suckthroughcc}
\end{lemma}

\begin{proof}
Consider the intersection of the sphere $E \cup D_i$ with the
projection plane.  By Lemma \ref{lemma:SintersectP}, we may assume
this intersection is a single simple closed curve $\gamma$, with one
arc of the curve running along the intersection of $D_i$ with the
projection plane, and the other arc meeting the knot strands exactly
in the intersections of $E$ with the knot strands.

If $C_i$ bounds a half--twist, we may isotope the diagram such that
all the crossings of the half--twist lie outside of the sphere $E\cup
D_i$, as in Figure \ref{fig:sphere-halftwist}.  In any case, the
portion of the diagram containing the sphere consists of $m_i$ strands
running parallel, embedded on the projection plane, entering $D_i$,
and $n_i<m_i$ parallel strands, embedded on the projection plane,
exiting $E$.

\begin{figure}
\begin{center}
\input{figures/sphere-halftwist.pstex_t}
\end{center}
\caption{}
\label{fig:sphere-halftwist}
\end{figure}

Isotope to move the ball bounded by $E\cup D_i$ to the opposite side
of $C_i$, as in Figure \ref{fig:shift-under}.  Notice that we have a
new diagram, with the same crossing circles as before, and the same
numbers of knot strands running through each crossing circle, except
there are now $n_i$ strands running through $C_i$.

\begin{figure}
\begin{center}
\input{figures/shift-under.pstex_t}
\end{center}
\caption{}
\label{fig:shift-under}
\end{figure}

If $C_i$ does not bound a half--twist, then this is a standard diagram
of a generalized augmented link, and the lemma is proved.  If $C_i$
does bound a half--twist, then this is no longer a standard diagram,
since the half--twist associated with $C_i$ no longer remains in a
neighborhood of $C_i$ after the isotopy.  In this case, we perform a
flype on the region of the diagram bounded by $\gamma$.  That is, we
rotate $180^\circ$ in the opposite direction of the half twist.  This
cancels the crossing of the ``old'' half--twist associated with $C_i$,
and adds a half--twist in the immediate neighborhood of $C_i$ after
the isotopy, as in Figure \ref{fig:flype}.

\begin{figure}
\begin{center}
\input{figures/flype.pstex_t} \\
\input{figures/flype2.pstex_t}
\end{center}
\caption{}
\label{fig:flype}
\end{figure}

Finally, notice that the flype does not add any crossings, aside from
those of the ``new'' half--twist adjacent to $C_i$.  For within the
region bounded by $\gamma$, knot strands on the projection plane are
taken back to the projection plane under the flype.  Crossing circles
are taken to crossing circles.  Half--twists rotate $180^\circ$ to
become identical half--twists.  Outside the region bounded by
$\gamma$, the diagram does not change at all, except to shift the
half--twist from one side of $\gamma$ to the other.  Thus in this case
we have established the lemma.
\end{proof}

\begin{lemma}
Given a standard diagram of a generalized augmented link $L$ in $S^3$,
there exists a new diagram that is also standard, but for which $D_i$
meets the knot strands in the minimal number of points.
\label{lemma:prop1}
\end{lemma}

\begin{proof}
If the original diagram does not satisfy this property, then combining
Lemmas \ref{lemma:meetnoDi} and \ref{lemma:suckthroughcc} we obtain a
new diagram for which the number of strands running through a single
crossing circle has been reduced.  Notice that there are only a finite
number of crossing circles, and a finite number of strands running
through each crossing circle.  Thus we need only repeat a finite
number of times, and we are left with a diagram for which $D_i$ meets
the knot strands in the minimal number of points.
\end{proof}

The next lemmas will give property (2) of Definition
\ref{def:reduced}.

\begin{lemma}
Let $L$ be a generalized augmented link in $S^3$ with standard
diagram.  Suppose there is an incompressible annulus $A$ embedded in
$L$ with boundary $C_i$ and $C_j$, for some $i\neq j$.  Then $A$ is
isotopic to an annulus which does not meet $D_i$ or $D_j$, and
intersects $\cup\, {\rm int}(D_k)$ in a (possibly empty) collection of
core curves of the annulus.
\label{lemma:annuli-through}
\end{lemma}

That is, any such annulus must either completely run through a twist
region, or completely miss the twist region.

\begin{proof}
First, apply Lemma \ref{lemma:prop1} to the diagram of $L$.  We may
assume $L$ has standard diagram for which each $D_i$ meets the knot
strands in the minimal number of points.  

Consider the intersection of $A$ with some $D_k$.  This is a
collection of closed curves.  Suppose one of them, say $\alpha$,
bounds a disk $D$ in $A$.  Then the disk obtained by taking $D_k$
outside $\alpha$ and $D$ inside $\alpha$ must meet the knot strands
the same number of times as $D_k$, since $D_k$ meets the knot strands
the minimal number of times.  Thus $\alpha$ must bound a disk in
$D_k$.  Hence we can isotope off, reducing the number of
intersections.

Now suppose $A$ intersects $D_i$.  Then by the preceding paragraph,
any intersection must be an essential curve $\alpha$ in $A$.  Thus we
may isotope $C_i$ along $A$ to this curve of intersection, and then
push off slightly, reducing the number of intersections of $A$ with
$D_i$.  Repeat, until there are no intersections with $D_i$.
Similarly for $D_j$.
\end{proof}

\begin{lemma}
Let $L$ be a generalized augmented link in $S^3$.  Suppose there is an
embedded annulus $A$ in $S^3-L$ with boundary $C_i$ and $C_j$, $i\neq
j$.  Then there is a generalized augmented link $L'$ in $S^3$ such
that first, $S^3-L$ is homeomorphic to $S^3-L'$, second, there is a
one-to-one correspondence between crossing circles and knot strands of
$L$ and $L'$, third, twisting disks of a standard diagram of $L$ meet
the knot strands the same number of times as the corresponding
twisting disks in a standard diagram of $L'$, but the crossing circle
corresponding to $C_i$ does \emph{not} meet the knot strands of $L'$
in a half--twist.
\label{lemma:redundant-twists}
\end{lemma}

\begin{proof}
Start with a standard diagram of $L$, and consider $C_i$ and $C_j$.
If one of $C_i$ and $C_j$ does not encircle a half--twist, then $L$,
possibly with $C_i$ and $C_j$ switched, satisfies the conclusions of
the lemma.

So suppose both $C_i$ and $C_j$ encircle half--twists.  Isotope $A$ as
in Lemma \ref{lemma:annuli-through}.  Then $A \cup D_i\cup D_j$ is a
sphere $S$ in $S^3$.  We may isotope $C_i$ and $C_j$ such that the
half--twists bounded by these crossing circles are both on the outside
of $S$.

Now perform a flype on the inside of $S$ in the direction opposite the
half--twist at $C_i$.  We wish to analyze what happens to the diagram
after the flype.  First, consider the portions of the diagram inside
$S$.  These are rotated $180^\circ$.  This rotation takes strands on
the projection plane back to the projection plane, and takes crossing
circles to crossing circles without affecting the number of strands
running through each crossing circle.  Finally, the rotation takes
half--twists to identical half--twists, in the same direction. 

Outside of $S$, the flype does not affect any strands of the diagram,
except that it removes the half--twist at $C_i$ and adds a new
half--twist at $C_j$.  This will either cancel with the half--twist
already at $C_j$, or it will add to it, putting a full--twist at
$C_j$.  In the case the flype leaves a full--twist at $C_j$, we
replace the flyped diagram with one in which the full twist at $C_j$
has been removed, replaced by parallel strands with no crossings at
$C_j$, as in Figure \ref{fig:annulus-flype}.  This link is the link
$L'$.  The complement of $L'$ is homeomorphic to the complement of
$L$.

\begin{figure}
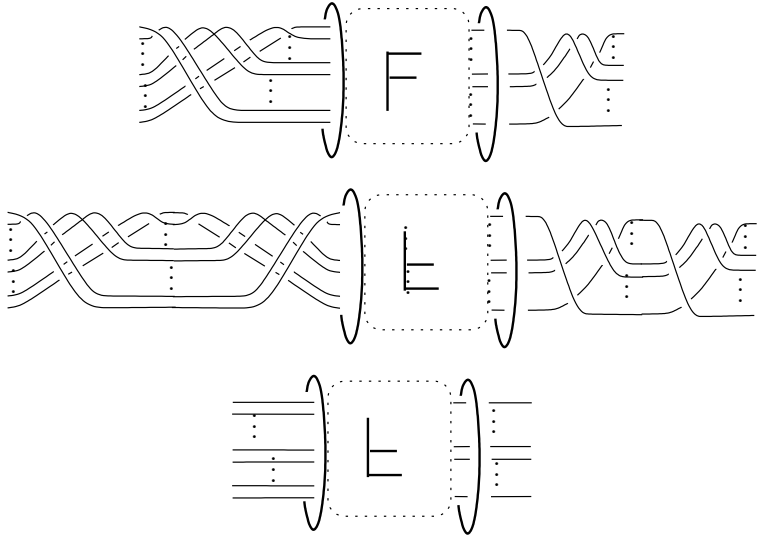

\input{figures/annulus-flype1.pstex_t}\\
\vspace{.1in}
\input{figures/annulus-flype2.pstex_t}\\
\vspace{.1in}
\input{figures/annulus-flype3.pstex_t}
\caption{Top to bottom: the link $L$ before the flype. After the
	flype. The link $L'$.}
\label{fig:annulus-flype}
\end{figure}

To finish the proof of the lemma, we need to show that the given
diagram of $L'$ is the standard diagram of a generalized augmented
link.

If no portion of the diagram of $L$ on the outside of $S$ crosses over
a portion of the diagram of $L$ on the inside of $S$, then this will
still be true for $L'$, and since both outside and inside of $S$ in
$L'$ we have portions of diagrams of generalized augmented links,
their union must be the diagram of a generalized augmented link.

So suppose some portion of the diagram of $L$ outside of $S$ crosses
over some portion diagram inside $S$.  These will form crossings of
the original standard diagram of $L$.  All crossings are associated
with half--twists.  Since the crossings of the half--twists associated
with $C_i$ and $C_j$ lie outside of $S$, and do not meet $S$, these
crossings must be associated to some $C_k$ with $k\neq i, j$.  Then we
may assume that $S$ intersects the corresponding twisting disk $D_k$,
or we may isotope $S$ out of a neighborhood of $D_i$ as in Lemma
\ref{lemma:SintersectP}, reducing the number of crossings of the outside of
$S$ with the inside of $S$.

So assume $S$ intersects $D_k$.  By Lemma \ref{lemma:annuli-through},
$D_k$ must meet the annulus $A$ in its non-trivial core curve.  Then
the portion of the diagram of $L$ inside $S$ at this intersection will
be a half--twist.  When we perform the flype, the half--twist becomes
an identical half--twist in the same direction.  Thus in a
neighborhood of $D_k$ inside of $S$, the diagram remains unchanged
after the flype.  Since the diagram does not change outside of $S$, a
neighborhood of $D_k$ will still be a half--twist after the flype.
Thus in all cases, the diagram of $L'$ is the standard diagram of a
generalized augmented link.
\end{proof}

\begin{proof}[Proof of Theorem \ref{thm:reduced-knot}]
Let $K$ be a knot in $S^3$ with a maximal twist region selection, and
let $L$ the corresponding generalized augmented link with standard
diagram.  If $L$ is not reduced, then it fails one of properties (1),
(2), or (3) of Definition \ref{def:reduced}.  Since our knot is
embedded in $S^3$, which has empty boundary, condition (3) will not
hold, so we need only show (1) and (2).

Suppose that property (1) fails.  Then by Lemma \ref{lemma:prop1}, we
may replace our diagram by a new standard diagram in which $D_i$
meets the knot strands in the minimal number of points.  Note that
\emph{a priori}, there may now be a twisting disk $D_i$ that meets the
knot strands in $m_i$ points, where $m_i \leq 1$.  In this case,
simply remove $C_i$ from the diagram, since twisting along $D_i$ when
$m_i \leq 1$ gives a manifold homeomorphic to that obtained by
performing meridional Dehn filling on $C_i$.

Suppose property (2) fails.  Then there is an annulus embedded in
$S^3-L$ with boundary on some $C_i$ and $C_j$.  By Lemma
\ref{lemma:redundant-twists}, we may replace the diagram of $L$ with a
standard diagram $D'$ of a generalized augmented link $L'$, without
increasing the numbers of crossing circles or strands running through
crossing circles, such that now the crossing circle corresponding to
$C_i$ doesn't encircle a half--twist in the diagram $D'$.  Consider
the image of the annulus in $S^3-L'$ under the homeomorphism of
$S^3-L$ and $S^3-L'$.  This is an annulus embedded in $S^3-L$ with one
boundary component on $C_i$.  We may isotope $C_i$ along this annulus
to be parallel to the other boundary component, $C_j$.  The result has
diagram $D''$ identical to $D'$, except that $C_i$ has been removed
and a crossing circle parallel to $C_j$ has been added.

Now, to consider the knot $K$, we perform twisting on $C_i$ and on
$C_j$.  However, note that both add full--twists to the same
generalized twist region.  Hence we may remove $C_i$ from the diagram,
and adjust the amount of twisting at $C_j$ to obtain the same link
$K$.

Perform Dehn filling along the crossing circles, adding the
appropriate number of twists so that the result has complement
homeomorphic to $S^3-K$.  This has a diagram given by replacing the
crossing circles $C_k$ in the diagram with a number of full--twists.
The diagram will be reduced, by construction.  Moreover, since knots
are determined by their complements \cite{gordon-luecke}, the diagram
is of a knot isotopic to $K$.
\end{proof}

\begin{remark}
Three remarks on the above proof.

First, note that when we reduced the diagram, we only removed crossing
circles from the unreduced augmented link, never added.  Thus the
resulting diagram will have at most as many generalized twist regions
as the original.  This is useful for Theorem
\ref{thm:gromov-norm-knots} below.

Second, what would be more useful would be if the reduced diagram had
at least as many half twists per generalized twist region as the
original, since many of our results in this paper require a knot with
a high number of half twists per twist region.  However, this may not
be the case if property (2) of Definition \ref{def:reduced} does not
hold for the original diagram.  In that case, we must concatenate two
generalized twist regions into a single one.  This may cancel half
twists in opposite directions, reducing the total number of half
twists in the new generalized twist region.

Finally, notice that the above proof used the fact that $K$ was a knot
only in the very last step.  To show the same result for links in
$S^3$, one needs to show that an isotopy of the generalized augmented
link followed by the appropriate Dehn filling is equivalent to an
isotopy of the original link.  We believe this is true (and possibly
even known), but we have not worked out the details here.
\end{remark}

\section{Applications to knots}
\label{sec:knot}
In this section, we apply the results of the previous sections to give
geometric information on knots in $S^3$.  Throughout, we will let $K$
be a knot in $S^3$, and let $D$ be a twist reduced diagram of $K$.  By
Theorem \ref{thm:reduced-knot}, we may assume any knot admits such a
diagram.  

\subsection{Torus decomposition and geometric type}
Given a knot $K$ in $S^3$, with twist reduced diagram $D$, form an
augmented link $L$ by adding crossing circles to the diagram in a
maximal twist region selection.  We have restrictions on essential
tori in $S^3-L$.

\begin{lemma}
  Let $T$ be an essential torus in $S^3-L$.  Then $T$ bounds a solid
  torus $V$ such that
  \begin{enumerate}
  \item $V$ is invariant under the involution $\sigma$,
  \item $V$ contains the link component $K$,
  \item $V\cap P$ is nonempty, with $P \cap \partial V$ containing no
    meridians of $V$.
  \end{enumerate}
\label{lemma:V}
\end{lemma}

\begin{proof}
Because $T$ is essential, it must intersect the surface $P$.  By the
invariant torus theorem \cite{holzmann}, we may assume $T$ is
invariant under $\sigma$.  By the solid torus theorem, $T$ bounds a
solid torus $V$ in $S^3$, which must also intersect $P$.  Because
$\sigma$ fixes $P\cap V$ pointwise, $\sigma$ must preserve $V$.  This
gives the first item of the lemma.

Suppose $V$ does not contain $K$.  At least one crossing circle must
be inside $V$, else a meridian of $V$ is a compressing disk of $T$,
contradicting incompressibility of $T$.  Let $C_i \subset V$.  The
circle $C_i$ bounds $D_i$.  If $D_i$ does not intersect $\partial V$,
then $D_i$ cannot intersect $K$, contradicting the definition of a
generalized augmented link.  Thus $D_i$ must intersect $\partial V$.

We may assume all such intersections are nontrivial curves on
$\partial V$, else we may replace disk portions of $D_i$ with disks of
$\partial V$.  This gives us a new $D_i$, still invariant under
$\sigma$ (since $V$ is invariant under $\sigma$), with required
properties of the definition of a generalized augmented link.  So we
assume all intersections of $D_i$ with $\D V$ are nontrivial.

Consider a curve of $D_i \cap \D V$ that is innermost in $D_i$.  This
bounds a disk $E$ on $D_i$.  The disk $E$ cannot lie in $V$, or it
would be a compressing disk.  Thus it must lie outside $V$.  But then
$V$ is unknotted in $S^3$.  Therefore replace $V$ with the exterior of
$V$ in $S^3$.  This is a solid torus that contains $K$.  Thus we have
the second item of the lemma.

Finally, suppose $\partial V \cap P$ consists of meridians of
$\partial V$, i.e., curves that bound disks in $V$.  Then since $V$ is
fixed by $\sigma$, $V \cap P$ must consist of two meridional disks.
The knot $K$ is contained in $V$, and $K$ is contained in $P$ by
Proposition \ref{prop:reflect}.  Since the $C_i$ are contained in a
neighborhood of $K$, the link $L$ must lie in a neighborhood of a
meridional disk of $V$.  But then $V$ is compressible.  Contradiction.
\end{proof}

Lemma \ref{lemma:V} allows us to classify Seifert fibered components
of the torus decomposition of $S^3-L$.  We have the following.

\begin{prop}
Let $K$ be a knot in $S^3$ with a twist reduced diagram.  Let $L$ be a
corresponding augmentation.  Let $\CT$ be the tori in the torus
decomposition of Theorem \ref{thm:tori}.  Then the tori of $\CT$ are
nested, bounding solid tori $V_1 \supset V_2 \supset \dots \supset V_k
\supset K$, and if $U$ is a component of $S^3-\CT$ which is not
homeomorphic to $T^2\times I$, such that $U-(L\cap U)$ is Seifert
fibered, then $U$ is the outermost component of the decomposition, and
$U-(L\cap U)$ is homeomorphic to the manifold of Figure
\ref{fig:mobius}.
\label{prop:sf-knot}
\end{prop}

\begin{proof}
Since each torus of $\CT$ is essential in $S^3-L$, by Lemma
\ref{lemma:V} each bounds a solid torus invariant under $\sigma$,
containing $K$, with $V\cap P$ a longitudinal slope.  Thus we can
arrange the solid tori in order of containment:
\[ K \subset V_k \subset \dots \subset V_1. \]

By Theorem \ref{thm:tori}, for each component $U$ of $S^3 - \CT$ which
is not homeomorphic to $T^2\times I$, $U-(L\cap U)$ is a reduced
augmented link complement.  Thus it has a boundary component
corresponding to some $C_j$.  Then any such component except the
outermost must have at least three boundary components: at least one
$C_j$, some $V_i$, and $V_{i+1}$ (or $K$).  Now, if this component is
Seifert fibered, Theorem \ref{thm:no-sf} tells us that $V_i$ and
$V_{i+1}$ must be parallel on $P$, and therefore $V_{i+1}$ cannot be a
subset of $V_i$.  This is impossible.

So only the outermost component may be Seifert fibered.  In this case,
again Theorem \ref{thm:no-sf} tells us there is just one $C_1$ in
the component $(S^3-L)-V_1$, and $(S^3-V_1) - C_1$ is a solid torus.
Since there is only one other link component, namely $V_1$, it must be
embedded in the M\"obius band $P$, parallel to the boundary of $P$, as
in Figure \ref{fig:mobius}.
\end{proof}

\begin{corollary}
Let $K$ be a torus knot with a twist reduced diagram $D$ and a maximal
twist region selection in which each twist region admits at least 6
half twists.  Then $K$ is a $(2, p)$ torus knot, and $D$ has one twist
region.
\label{cor:torus-knot}
\end{corollary}

\begin{proof}
$S^3-K$ is Seifert fibered.  Since $D$ is twist reduced, adding
crossing circles to twists of $D$ yields a reduced augmented link $L$
in $S^3$.  By Theorems \ref{thm:tori} and \ref{thm:dehnfill}, there is
a sublink $\hat{L}$ of $L$, possibly consisting of fewer crossing
circles, and a collection of tori $\hat{\TT}$ such that that
components of $(S^3-\hat{L}) -\hat{\TT}$ are reduced augmented links
with the same geometric type as those of $S^3-K$.  Thus the components
must all be Seifert fibered.

By Proposition \ref{prop:sf-knot}, only the outermost component is
Seifert fibered, and it is of the form of Figure \ref{fig:mobius}.
Thus the collection of tori of $\hat{\TT}$ is empty.  When we twist to
insert at least six half twists, we obtain a $(2,p)$ torus knot.
\end{proof}

\begin{theorem}
Let $K$ be a knot in $S^3$ which is toroidal, with a twist--reduced
diagram and a maximal twist region selection with at least $6$
half--twists in each generalized twist region.  Let $L$ denote the
corresponding augmentation.  Then there exists a sublink $\hat{L}$ of
$L$, possibly containing fewer crossing circles, such that:
\begin{enumerate}
	\item The essential tori of the JSJ decomposition of $S^3-K$ are in
	one--to--one correspondence with those of $S^3-\hat{L}$.
	\item Corresponding components of the torus decompositions have the
	same geometric type, i.e. are hyperbolic or Seifert fibered.
	\item The essential tori of $S^3-\hat{L}$ and $S^3-K$ form a
	collection of nested tori, each bounding a solid torus in $S^3$
	which contains $K$, and which is fixed under the reflection of
	$S^3-L$.
\end{enumerate}
\label{thm:toroidal}
\end{theorem}

\begin{proof}
By Theorems \ref{thm:tori} and \ref{thm:dehnfill}, there is some
sublink $\hat{L}$ of $L$, possibly containing fewer crossing circles
(i.e. those inside any $T^2\times I$ components), and tori $\hat{\TT}$
such that the tori form a torus decomposition of $S^3-K$, and
corresponding components of the decompositions of $S^3-K$ and
$S^3-\hat{L}$ have the same geometric type.  By Lemma \ref{lemma:V},
item (3) must hold for all essential tori.

All that remains to prove is that the decomposition given by
$\hat{\mathbb T}$ is the JSJ decomposition, i.e. it is the unique
minimal torus decomposition of $S^3-\hat{L}$ and $S^3-K$.  If not,
then some torus $T$ of $\hat{\mathbb T}$ separates two Seifert fibered
components.  But by Proposition \ref{prop:sf-knot}, this is
impossible: only the outermost component of $S^3-\hat{L}$ can be
Seifert fibered.  Thus any other components must be hyperbolic, and so
$\hat{\mathbb T}$ is the unique torus decomposition of $S^3-\hat{L}$.

Similarly, the only way $\hat{\mathbb T}$ can fail to be the unique
torus decomposition of $S^3-K$ is if some essential torus $T$ of
$\hat{\mathbb T}$ splits a single Seifert fibered component into two.
But again the components of $S^3-K$ have the same geometric type as
those of $S^3-\hat{L}$, by Theorem \ref{thm:dehnfill}.  Hence
$\hat{\mathbb T}$ must be the unique torus decomposition of $S^3-K$.
\end{proof}

\section{Application: Gromov norms}
In this section, we apply the previous results to bound the Gromov
norms of many toroidal knots.  The results use heavily the particular
torus decompositions developed in previous sections, as well as
results on hyperbolic generalized augmented links in
\cite{purcell:aug}.

First, we insert a bound on the volume of a hyperbolic augmented link
in a solid torus, whose proof follows immediately from
\cite{purcell:aug}.

\begin{prop}
Let $L$ be a generalized augmented link in a solid torus $V$, with
$t_i$ crossing circles.  Suppose $V-L$ is hyperbolic.  Then its volume
satisfies $\vol(V-L) \geq 2 \, v_8 \, t_i,$ where $v_8 \approx
3.66386$ is the volume of a hyperbolic regular ideal octahedron.
\label{prop:vol-solid-torus}
\end{prop}

\begin{proof}
By Lemma 4.1 of \cite{purcell:aug}, the volume of a hyperbolic
generalized augmented link in $S^3$ with $t$ generalized twist regions
is at least $2\,v_8(t-1)$.  We may use this result to bound volumes of
hyperbolic generalized augmented links in a solid torus as follows.
Embed the solid torus as one of the solid tori of a standard genus one
Heegaard splitting of $S^3$.  Let $C$ be the core of the other solid
torus.  Then the embedding gives a generalized augmented link in $S^3$
where the curve $C$ is an additional crossing circle component.  Thus
the volume of a hyperbolic augmented link in a solid torus is at least
$2\,v_8 (t_i+1 -1) = 2\,v_8\,t_i$, where $t_i$ denotes the number of
crossing circles.
\end{proof}

Alternately, Proposition \ref{prop:vol-solid-torus} may be proved by
cutting $V-L$ along $P$ and crossing disks, and using this to
determine the number of vertices and edges of a one--skeleton as in
the proof of Lemma 4.1 of \cite{purcell:aug}.

\begin{theorem}
Let $K$ be a knot with twist reduced diagram with $t$ generalized
twist regions.  Let $L$ be the corresponding generalized augmented
link and let $t_0$ denote the number of components of the torus
decomposition of the form $(T^2\times I) - C_i$.  Then the Gromov norm
of $S^3-L$ satisfies
$$\| [S^3-L] \| \geq 2\,v_3\,v_8 (t - t_0 -1),$$ where $v_3 \approx
1.0149$ is the volume of a regular ideal hyperbolic tetrahedron,
and $v_8 \approx 3.66386$ is the volume of a regular ideal
hyperbolic octahedron.
\label{thm:gromov-norm-links}
\end{theorem}

\begin{proof}
By Proposition \ref{prop:sf-knot} and Theorem \ref{thm:tori}, $S^3-L$
admits a torus decomposition such that each piece of the decomposition
is either a link in $T^2 \times I$, or a generalized augmented link in
a solid torus.  Moreover, the only Seifert fibered pieces each meet
just one crossing circle component, either as the outermost component
of the decomposition or as a link in $T^2\times I$.  Since the Gromov
norm is $v_3$ times the sum of the volumes of the hyperbolic pieces
\cite{soma}, we use the results of Proposition
\ref{prop:vol-solid-torus}.  

Now, in case the outermost piece is hyperbolic, it is an augmented
link in $S^3$ and so its volume is at least $2\,v_8 (t_1-1)$, where
$t_1$ denotes the number of crossing circles in the outermost piece.
Then the Gromov norm is at least
$$\| [S^3-L] \| \geq v_3\left( 2\,v_8(t_1 -1) + \sum
2\,v_8\,t_i\right) = 2\,v_3\,v_8\,(t - t_0 -1),$$ where the sum in the
center is over all components of the torus decomposition which are
homeomorphic to hyperbolic generalized augmented links in a solid
torus.

If the outermost piece is hyperbolic, then it contributes nothing to
the Gromov norm.  By Proposition \ref{prop:sf-knot}, it contains just
one crossing circle, and so the Gromov norm is at least
$$\| [S^3-L]\| \geq v_3\,\sum 2\,v_8\,t_i = 2\,v_3\,v_8\,(t - t_0 -1),$$
where the final $-1$ corresponds to the single crossing circle in the
outermost component of the torus decomposition.
\end{proof}

For knots:

\begin{theorem}
Let $K$ be a knot in $S^3$ which is toroidal, with a twist--reduced
diagram at least $7$ half--twists in each generalized twist region.
Let $L$ denote the corresponding augmentation, and let $\hat{L}$
denote the sublink of Theorem \ref{thm:toroidal}.  Let $t$ denote the
number of crossing circles of $\hat{L}$.  Then the Gromov norm of
$S^3-K$ satisfies
$$\| [S^3-K] \| \geq 0.65721\,(t-1).$$
\label{thm:gromov-norm-knots}
\end{theorem}

\begin{proof}
Since the Gromov norm is $v_3$ times the sum of volumes of hyperbolic
pieces of the torus decomposition, we bound volumes of the hyperbolic
pieces.  Each is obtained by Dehn filling a generalized augmented
link.  Moreover, by \cite[Proposition 3.5]{purcell:aug}, the length of
the slope of Dehn filling is at least $\sqrt{(1/4) + c_i^2}$, where
$c_i$ is the number of half--twists.  Thus if there are at least $7$
half--twists in each twist region, then the Dehn filling slopes have
length at least $\sqrt{49.25} \approx 7.0178 > 2\pi$.

Now apply \cite[Theorem 1.1]{fkp}.  This theorem bounds the volume
under Dehn filling in terms of the length of the filling slope and the
volume of the unfilled manifold.  In particular, for twisting a
hyperbolic augmented link with $t_i$ crossing circles in a solid
torus, by Proposition \ref{prop:vol-solid-torus} we have
$$\mbox{volume after twisting} \geq \left(1 -
\left(\frac{2\pi}{\sqrt{49.25}}\right)^2\right)^{3/2} \,(2\,v_8\,t_i)
> 0.64756\,t_i.$$

Only the outermost component of the torus decomposition may be Seifert
fibered.  If it is not, the volume of the outermost piece is at least
$0.64756\,(t_1-1)$, where $t_1$ is the number of crossing circles, by
\cite[Theorem 4.2]{purcell:aug}.  Since $\hat{L}$ is obtained from $L$
by removing all crossing circles which lie in components $(T^2\times
I)-C_j$ of the torus decomposition of $S^3-L$ (Theorem
\ref{thm:tori}), the total number of crossing circles in $\hat{L}$ is
$t_1 + \sum t_i$, where the sum is over all pieces homeomorphic to
hyperbolic generalized augmented links in solid tori.  Thus the Gromov
norm satisfies
$$\|[S^3-K]\| \geq v_3 \left( 0.64756\,(t_1-1) + \sum 0.64756\,t_i
\right) = 0.65721\ldots\,(t-1).$$

If the outer most component is Seifert fibered, then the outermost
component contains a single crossing circle component, by
Proposition \ref{prop:sf-knot}, and so the Gromov norm satisfies
$$\|[S^3-K]\| \geq v_3 \, \sum 0.64756\,t_i = 0.65721\ldots (t-1).$$
\end{proof}


\bibliography{references}

\end{document}